\documentclass[twoside]{article}
\usepackage{amsmath,amssymb,amsthm} 
\usepackage{rotating}
\usepackage{adjustbox}
\usepackage{array}      
\usepackage{multicol}         
\usepackage[bottom]{footmisc}  
\usepackage{float}
\usepackage{cancel}
\usepackage[usenames,dvipsnames]{pstricks}
\usepackage{pst-grad} 
\usepackage{pst-plot}  
\usepackage{multirow} 

\usepackage{color}
\usepackage{setspace}
\usepackage{colortbl}
\usepackage{wrapfig}
\usepackage{lmodern}
\usepackage{tcolorbox}
\usepackage[pangram]{blindtext}
\usepackage{graphicx,epstopdf}

\usepackage{fancyhdr}
\usepackage{parskip}
\usepackage{colortbl,longtable}
\usepackage[explicit,calcwidth]{titlesec}
\usepackage[left=4.1cm, right=4.1cm, top=5.0cm, bottom=4.50cm]{geometry}
\usepackage{comment}
\usepackage{breakcites}

\def\R{\mathbb{R}}

\newcommand{\f}[2]{\frac{\displaystyle{#1}}{\displaystyle{#2}}}

\titleformat{\section}
  {\addtolength{\titlewidth}{2pc}\normalfont\bfseries}
  {\thesection.\,}
  {0.19em}{#1} 
  []
  
\titleformat{\subsection}
  {\addtolength{\titlewidth}{2pc}\normalfont\bfseries}
  {\thesubsection.\,}
  {0.19em}{#1} 
  []
  
\newtheorem{theorem}{Theorem}[section]

\newtheorem{corollary}{Corollary}[section]
\newtheorem{remark}{Remark}[section]
\newtheorem{lemma}{Lemma}[section]

\begin{document}

\title{}
\date{}
\maketitle
\vspace{-4.5cm}

\begin{center}
{\bf \large SIR models with vital dynamics, reinfection and\\[0.1cm] randomness to investigate the spread of infectious diseases}
\end{center}
\vspace{0.25cm}
\begin{center}
Javier L\'opez-de-la-Cruz$^\text{\footnotesize{a,}}$\footnote[1]{corresponding author: javier.lopez.delacruz@upm.es} and Alexandre N. Oliveira-Sousa$^\text{\footnotesize{b}}$
\end{center}
\vspace{0.5cm}

\begin{minipage}{0.97\textwidth}
\begin{itemize}
\item[$^\text{\footnotesize{a}}$] {\small Dpto. de Matem\'atica Aplicada a las TIC, Escuela T\'ecnica Superior de Ingenieros Inform\'aticos, Campus de Montegancedo, Universidad Polit\'ecnica de Madrid, 28660 Boadi\-lla del Monte, Madrid, Spain.}
\item[$^\text{\footnotesize{b}}$] {\small Dpto. de Matem\'atica, Universidade Federal de Santa Catarina, 88040-900 Florian\'opolis, Santa Catarina, Brazil}
\end{itemize}
\end{minipage}

\vspace{0.7cm}

\begin{center}
\begin{minipage}{0.97\textwidth}
\begin{abstract}
\small We investigate SIR models with vital dynamics, reinfection, and randomness at the transmission coefficient and recruitment rate. Initially, we conduct an extensive analysis of the autonomous scenario, covering aspects such as local and global well-posedness, the existence and internal structure of attractors, and the presence of gradient dynamics. Subsequently, we explore the implications of small nonautonomous random perturbations, establishing the continuity of attractors and ensuring their topological structural stability. Additionally, we study scenarios in which both the transmission coefficient and the recruitment rate exhibit time-dependent or random behavior. For each scenario, we establish the existence of attractors and delineate conditions that determine whether the disease is eradicated or reaches an endemic state. Finally, we depict numerical simulations to illustrate the theoretical results. 

\end{abstract}
\vspace{0.75cm}

\end{minipage}
\end{center}
\vspace{0.5cm}








\section{Introduction}\label{introduction}

Epidemics have influenced the humanity over the history, bringing about consequences as catastrophic as the ones caused by wars and making entire populations to become extinct. 
More recent epidemics are the Avian influenza, 
the Hong Kong flu, 
and the COVID-19, that we have been suffering since the beginning of 2020, due to which more than 6.9 million people have lost their lives in the World up to now.
Apart from the human losses, 
epidemics also have important social consequences and have a huge negative impact on the economy, see \cite{miranda2022}.  
In addition, many infectious diseases, which have caused epidemics over history, are not eradicated yet and new ones are arising over the years. Hence, it is essential to set up mathematical models to describe the evolution of the diseases and then understand their evolution and make decisions to control the spread of the disease as soon as possible minimizing the consequences.

The first paper known about mathematical modelling of epidemics was published in 1760 (see \cite{daniel-bernoulli-1960}), where D. Bernoulli studied the spread of the Smallpox by means of system of ordinary differential equations. However, the study of mathematical models in Epidemiology was not developed until the very beginning of the twentieth century, when W. H. Hamer formulated in 1906 a discrete mathematical model to describe the transmission of the Measles (see \cite{william-heaton-hamer-1906}). Afterwards, in 1911, R. Ross provided a mathematical model to predict the spread of an outbreak of Malaria (see \cite{ronald-ross-1911}), where he proved that it was enough to reduce partially the population of mosquitoes to eradicate the epidemics.

Nevertheless, in 1927 Kermarck and McKendrick established a mathematical model to describe the spread of an infectious disease (see
\cite{kermack-mckendrick-1927})
which is the base in mathematical epidemiology ever since. This classical model, called SIR, is a compartmental model, i.e., it divides the population in three groups. The first group contains the susceptible individuals ($S$), people who are not infected but can contract the disease when contacting with an infected. The second group contains the infected individuals ($I$). The third group contains the recovered individuals ($R$), people who are recovered with immunity or who pass away because of the disease.

\begin{figure}[H]
\begin{center}
\includegraphics[width=0.3\textwidth]{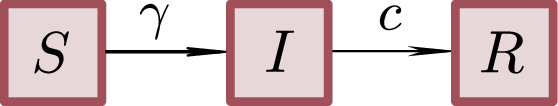}
\caption{Diagram of the classical SIR model}
\label{diagramaKermackMcKendrick}
\end{center}
\end{figure}

In Figure \ref{diagramaKermackMcKendrick} we depict a diagram with the flow of states regarding the classical SIR model, where $\gamma>0$ represents the transmission coefficient of the disease and $1/c>0$ describes the mean time that an individual remains infected when contracting the disease. As we can observe, the susceptible individuals only abandon their group when contracting the infection (after having contact with a infected individual) and the infected ones can leave their groups when they are recovered or pass away due to the disease.

The classical SIR model by Kermack and McKendrick is given by the following differential system
\begin{align}
S^\prime&=-\gamma\frac{SI}{S+I+R},\label{Soriginal}\\
I^\prime&=\gamma\frac{SI}{S+I+R}-cI,\label{Ioriginal}\\
R^\prime&=cI,\label{Roriginal}
\end{align}
where $S=S(t)$, $I=I(t)$ and $R=R(t)$ denote the number of susceptible, infected and recovered individuals at time $t$, respectively. 

Unfortunately, the classical SIR model \eqref{Soriginal}-\eqref{Roriginal} assumes strong restrictions. For instance, it assumes that the total population $N=S+I+R$ is constant. However, infectious diseases do not disappear soon usually (as we see with COVID-19), then the size of the population usually varies due to external factors. Another important restriction is to assume that every infected individual either recovers with immunity or passes away (this is not the case with COVID-19).

\begin{figure}[H]
\begin{center}
\includegraphics[width=0.3\textwidth]{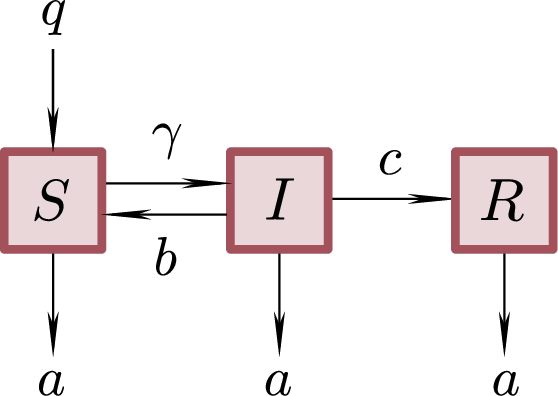}
\caption{Diagram of the modified SIR model}
\label{diagram}
\end{center}
\end{figure}

To set up a more realistic model and avoid the drawbacks pointed out above, we consider vital dynamics, i.e., we introduce the recruitment rate $q>0$, which takes into account the new susceptible individuals arising in the population during the epidemic (such as births) and another parameter, $a>0$, to take also into account the deaths that are not caused by the disease (such as natural deaths, accidents or others). 
In addition, we introduce another parameter $b>0$ to describe those infected individuals that recover without immunity. The diagram with the flow of states for the modified SIR model is in Figure \ref{diagram}.

The differential system corresponding to the modified SIR model is given by
\begin{align}
S^\prime&=q-aS+bI-\gamma\frac{SI}{S+I+R},\label{Snew}\\
I^\prime&=-(a+b+c)I+\gamma\frac{SI}{S+I+R},\label{Inew}\\
R^\prime&=cI-aR.\label{Rnew}
\end{align}

In \cite{Kloeden-Potzsche} the authors study some two-dimensional versions of system \eqref{Snew}-\eqref{Rnew} with a time-dependent variable population $q=q(t)$ and 
they motivate the case where $\gamma=\gamma(t)$ is also time-dependent and satisfies a mean zero condition in time, a condition which we will use for the random case. 
More precisely, they study the disease-free ``non-autonomous equilibrium'' of the systems. Usually, these nonautonomous equilibria are special global solutions that resemble the behavior of the autonomous equilibria in nonautonomous problems, see  \cite{Carvalho-Langa-2,Chueshov}. Authors typically analyze the stability of the disease-free equilibrium when dealing with SIR models: if it is stable, the disease is eradicated; however, if it is unstable, the disease becomes endemic, see also \cite{Kloeden-Kozyakin}. Here, we also want to study if this solution is hyperbolic in the sense that the linearized system admits an exponential dichotomy, see \cite{Caraballo-Carvalho-Langa-OliveiraSousa-2021,Carvalho-Langa-2}.

On the other hand, the authors in \cite{Caraballo-Colucci-17} focus on the system \eqref{Snew}-\eqref{Rnew} in a random setting. First, $q$ is considered as a bounded random variable. After proving the existence of random attractors, they provide conditions under which the disease is eradicated or becomes endemic. Moreover, they study the model under environmental stochastic perturbations and prove the existence of random attractors. 

Inspired by those works, we introduce the following modifications to set up even more realistic systems. First, we propose to replace $\gamma$ with $\gamma+\Phi(\theta_t\omega)$, where $\Phi$ is a random variable driven by a random flow $\theta$. This idea is further supported by \cite{Albani-Zubelli}, elucidating that the transmission coefficient demonstrates randomness in practical scenarios.
The resulting random model is given by
\begin{align}
S^\prime&=q-aS+bI-(\gamma+\Phi(\theta_t\omega))\frac{SI}{S+I+R},\label{S_random_gamma}\\
I^\prime&=-(a+b+c)I+(\gamma+\Phi(\theta_t\omega))\frac{SI}{S+I+R},\label{I_random_gamma}\\
R^\prime&=cI-aR\label{R_random_gamma}.
\end{align}

After that, as a second step, we propose to consider the recruitment rate $q$ in system \eqref{S_random_gamma}-\eqref{R_random_gamma} as a nonautonomous and random variable $q(t,\theta_t\omega)$, based on the stochastic nature of new individuals joining the population (see \cite{Caraballo-Colucci-17}). 
We investigate the dynamics of the disease in both systems  
by employing the theory of nonautonomous and random dynamical systems. We aim to find conditions under which the disease is eradicated or becomes endemic, the existence of attractors, and hyperbolic solutions. 





This paper is organized as follows. In Section \ref{SIRautonomous}, we investigate the SIR autonomous model with vital dynamics and reinfection \eqref{Snew}-\eqref{Rnew} proving that it is associated with a gradient semigroup. 
As an application, in Section \ref{SIR-small-nonautrandom-perturbations}, we study small non-autonomous random perturbations of \eqref{Snew}-\eqref{Rnew} and prove a result of continuity and topological structural stability applying results from 
\cite{Caraballo-Carvalho-Oliveira-Sousa-NRA}. 
In Section \ref{SIRgamma}, we study the system \eqref{S_random_gamma}-\eqref{R_random_gamma}, we obtain attractors existence and conditions under which the disease is eradicated. In Section \ref{SIRqgamma}, we consider the random model \eqref{S_random_gamma}-\eqref{R_random_gamma} with a nonautonomous random $q(t,\omega)$ providing similar results to the ones in the previous section, and we study hyperbolicity for the disease-free solution. Finally, in Section \ref{sec:numerical_simulations} we give examples of specific noises that can be used and illustrate the theoretical work. 

\section{Autonomous SIR model with vital dynamics and reinfection}\label{SIRautonomous}

In this section, we investigate the asymptotic behavior of \eqref{Snew}-\eqref{Rnew}. We prove the local and global well-posedness, such that the solutions generate a semigroup that possesses a global attractor. We also study the internal structure of the attractors, and prove that the semigroup is Gradient, so there is a Lyapunov function which is non-increasing along orbits and constant in each equilibrium, see \cite{Aragao-Costa-2011}. 
Before we proceed, we briefly recall some basic notions of autonomous dynamical systems, and we recommend \cite{Bortolan-Carvalho-Langa-book,Carvalho-Langa-Robison-book} for an introduction to this topic. 

A \textbf{semigroup} $\{T(t): t\geq 0\}$ in a metric space $(X,d)$ is a family of continuous operators such that $T(0)=Id_X$, $T(t)T(s)=T(t+s)$ for $t,s\geq 0$, and the mapping
$(t,x)\mapsto T(t)x$ is continuous. 
We say that $\xi:\mathbb{R}\to X$ is a \textbf{global    solution} for $\{T(t): t\geq 0\}$ if $T(t)\xi(s)=\xi(t+s)$ for all $t\geq 0$, $s\in \mathbb{R}$. Finally, a \textbf{global attractor} for $\{T(t): t\geq 0\}$ is a compact subset $\mathcal{A}$ of $X$ such that $\mathcal{A}$ is invariant ($T(t)\mathcal{A}=\mathcal{A}$, for all $t\geq 0$) and attracts every bounded subset of $X$, i.e., for any bounded subset $B\subset X$, 
\begin{equation*}
    \lim_{t\to +\infty} dist(T(t)B,\mathcal{A})=0,
\end{equation*}
\noindent where $\displaystyle dist(A,B)=\sup_{a\in A}d(a,B)$ is the Hausdorff semi-distance between $A$ and $B$.

Now, we study the set of equilibria of system \eqref{Snew}-\eqref{Rnew}. If $\gamma\leq a+b+c$, it has only the {\it disease-free equilibrium} $E_0^*=(q/a,0,0)$. However, the so-called \textit{endemic equilibrium} $E_1^*=(S^*,I^*,R^*)$ appears in the case $\gamma>a+b+c$, where
\begin{equation}\label{eq-equilibrium_E_1^*}
S^*=\frac{q(a+b+c)}{\gamma a},\quad 
I^*=\frac{q(\gamma -a-b-c)}{\gamma(a+c)},\quad
R^*=\frac{c}{a}I^*.
\end{equation}

In the following, we study hyperbolicity for the equilibria, using standard results from the theory of ODEs, see for instance \cite{Martcheva}.
\
\begin{lemma}\label{lemma-E_0-E_1-hyperbolics}

Let $E_0^*$ and $E_1^*$ be defined as above.
\begin{enumerate} 
    \item If $\gamma\leq a+b+c$, then $E_0^*$ is the only equilibirum of system \eqref{Snew}-\eqref{Rnew}.
    Moreover, if $\gamma< a+b+c$, then $E_0^*$ is exponentially stable.
    \item If $\gamma>a+b+c$, then $E_0^*$ is a sadlle point with one-dimensional local unstable manifold, and  $E_1^*$ is a exponencially stable equilibrium.
\end{enumerate} 
\end{lemma}

\begin{proof} 
It is straightforward to verify that  the Jacobian matrix of the vector field of system \eqref{Snew}-\eqref{Rnew} is given by
\begin{equation}\label{jacobianmatrix}
\left(\begin{array}{ccc}
\displaystyle{-a-\gamma\frac{I(I+R)}{(S+I+R)^2}} & 
\displaystyle{b-\gamma\frac{S(S+R)}{(S+I+R)^2}} &
\displaystyle{\gamma\frac{SI}{(S+I+R)^2}} \\[2.5ex]
\displaystyle{\gamma\frac{I(I+R)}{(S+I+R)^2}} & \displaystyle{-(a+b+c)+\gamma\frac{S(S+R)}{(S+I+R)^2}}
 & \displaystyle{-\gamma\frac{SI}{(S+I+R)^2}}
 \\[2.5ex]
0 & c & -a
\end{array}\right).
\end{equation}
On the one hand, if $\gamma\leq a+b+c$, the eigenvalues of the Jacobian matrix \eqref{jacobianmatrix} at $E^*_0$ are given by $\lambda_1=-a<0$,  with algebraic and geometric multiplicity $2$, and $\lambda_2=\gamma-a-b-c$, with algebraic and geometric multiplicity $1$. Hence, the unique equilibrium $E_0^*$ is exponentially stable when $\gamma<a+b+c$ and it is no longer hyperbolic if $\gamma =a+b+c$.

%

 On the other hand, if $\gamma>a+b+c$, the Jacobian matrix \eqref{jacobianmatrix} at $E_0^*$ has one positive eigenvalue 
$\lambda_2=\gamma-a-b-c>0$, then $E_0^*$ becomes a saddle-point. Finally, it is straightforward to verify, thanks to the Routh–Hurwitz Criteria (see \cite[Theorem 5.1]{Martcheva}), that the roots of the characteristic polynomial associated to the Jacobian matrix \eqref{jacobianmatrix} at $E_1^*$ are negative or have negative real part, whence $E_1^*$ is exponentially stable.
\end{proof}  

We denote $\R^3_+=\{(S,I,R)\in \mathbb{R}^3:  S,I,R\geq 0\}$. The next result characterizes the dynamics of the autonomous system in detail.  

\begin{theorem}\label{th-SIR-autonomouos}
For any $u_0=(S_0,I_0,R_0)\in\R^3_+$, system \eqref{Snew}-\eqref{Rnew} has a unique global solution $u(\cdot,u_0)=(S(\cdot, u_0), I(\cdot,u_0), R(\cdot,u_0))\in  {C}^1([0,+\infty);\R^3_+)$, which generates a semigroup $\{T(t): t\geq 0\}$ in $\R^3_+$, defined as $T(t)u_0=u(t,u_0)$, that posseses a global attractor $\mathcal{A} \subset \mathbb{R}^3_+$ contained in
\begin{equation}
B_0:=\left\{(S,I,R)\in\R^3_+\,:\, S+I+R= \f{q}{a} \right\}. \label{eq-B0qconstante-autonomo}
\end{equation}
Moreover, the following statements hold true:
\begin{enumerate}
    \item If $\gamma<a+b+c$, the global attractor is $\mathcal{A}=\{E_0^*\}$.
    \item If $\gamma>a+b+c$,
    the global attractor is
    $\mathcal{A}=\{E_0^*,E_1^*\}\cup \{\xi(t): t\in \mathbb{R}\}$, where $\xi$ is a global solution such that
    $E_0^* \stackrel{t\to-\infty}{\longleftarrow} \xi(t)\stackrel{t\to+\infty}{\longrightarrow} E_1^*$ (see Figure \ref{diagrama1}). In particular, the semigroup $\{T(t): t\geq 0\}$ is gradient with respect to the set of equilibria $\{E_0^*, E_1^*\}$.
    \begin{figure}[H]
    \begin{center}
    \includegraphics[width=0.47\textwidth]{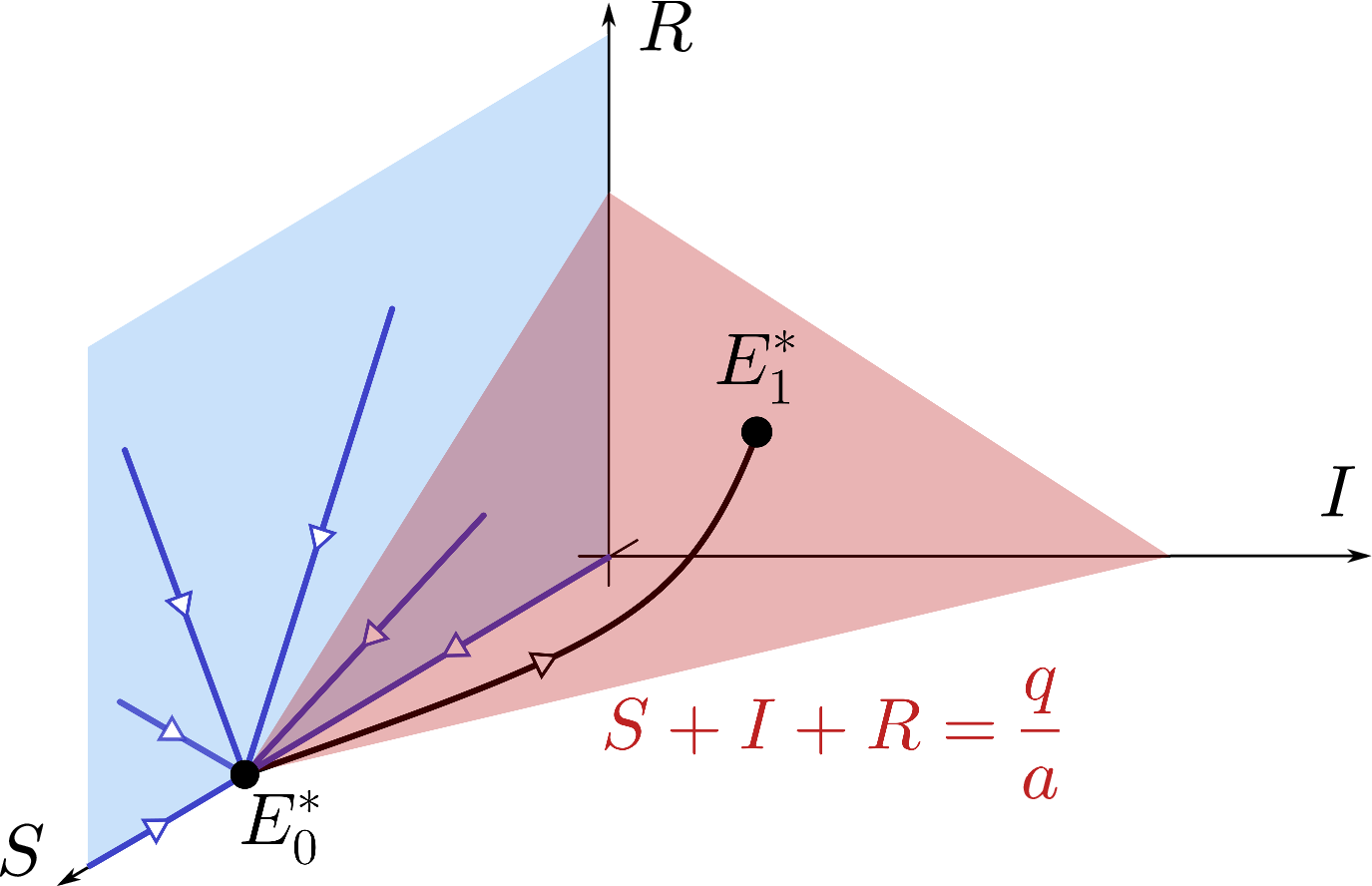}
    \end{center}
    \caption{Diagram of the autonomous case when $\gamma>a+b+c$.}
    \label{diagrama1}
\end{figure}
\end{enumerate}
\end{theorem}

\begin{proof}
The proof that the semigroup is well defined, the existence of the global attractor $\mathcal{A}$ and that it is contained in $B_0$ is a particular case of Theorem \ref{th-SIR-perturbing-gamma} in Section \ref{SIRgamma}.
From Lemma \ref{lemma-E_0-E_1-hyperbolics}, it is straightforward to verify that the global attractor is $\mathcal{A}=\{E_0^*\}$ and it attracts bounded sets with exponential rate if $\gamma<a+b+c$.

The proof that the $\omega$-limits and $\alpha$-limits consist only of equilibria is a consequence of the Poincar\'e–Bendixson Trichotomy and the Dulac-Bendixson Criterion (see \cite[Theorem 3.5 and Theorem 3.6]{Martcheva}). To see this, we use that $B_0$ is positively invariant under $\{T(t): t\geq 0\}$, then the system \eqref{Snew}-\eqref{Rnew} can be reduced to another one with only variables $S$ and $I$. Hence, it is enough to consider the Dulac function 
$D(S,I)=1/I$ to conclude the non-existence of non-constant periodic orbits.

We consider now $\gamma>a+b+c$. From Lemma \ref{lemma-E_0-E_1-hyperbolics}, the set of equilibria is $\{E_0^*,E_1^*\}$, where $E_1^*$ is exponentially stable and $E_0^*$ is a hyperbolic equilibrium with 
one-dimensional local unstable manifold. 
Thus there exists a unique global solution $\xi:\mathbb{R}\to \mathbb{R}^3_+$ such that $E_0^* \stackrel{t\to-\infty}{\longleftarrow} \xi(t)\stackrel{t\to+\infty}{\longrightarrow} E_1^*$. Therefore $\{E_0^*,E_1^*\}\cup \{\xi(t):t\in \mathbb{R}\} \subset \mathcal{A}$. Now let us prove that $\{E_0^*,E_1^*\}\cup \{\xi(t):t\in \mathbb{R}\}$ in fact the entire global attractor $\mathcal{A}$.

Recall that the global attractor $\mathcal{A}$ is the union of all bounded global solutions (see \cite[Theorem 1.7]{Carvalho-Langa-Robison-book}). Let $\bar{\xi}$ be a bounded global solution, then $\omega(\bar{\xi})\in\{E_0^*, E_1^*\}$. 
First, assume that $\omega(\bar{\xi})=E_0^*$. In this case $\bar{\xi}(\tau)\in W^s(E_0^*)$ for any $\tau \in \mathbb{R}$, where $W^s(E_0^*)$ denotes the stable manifold of $E_0^*$. However, $W^s(E_0^*)=\{(S,I,R)\in\R^3_+\,:\,I=0\}$. Since $\bar{\xi}(\mathbb{R})\subset \mathcal{A}$ and $\mathcal{A}\subset B_0$, then $\bar{\xi}(\mathbb{R})\subset B_0\cap W^s(E_0^*)=\{E_0^*\}$, thus $\bar{\xi}$ is constant equal to $E_0^*$. 
Now suppose that $\omega(\bar{\xi})=E_1^*$. Since $\alpha(\bar{\xi})\in\{E_0^*,E_1^*\}$ and $E_1^*$ is stable, we have $\alpha(\bar{\xi})=E_0^*$. Nevertheless, since the local unstable manifold of $E_0^*$ is one dimensional, the global solution $\bar{\xi}$ must coincide with $\xi$.

Hence $\mathcal{A}=\{E_0^*,E_1^*\}\cup \{\xi(t):t\in \mathbb{R}\}$ and by the same reasoning, we conclude that $\{T(t):t \in \R^+\}$ is dynamically gradient with respect to $\{E_0^*,E_1^*\}$, which is equivalent to say that $\{T(t):t \geq 0\}$ is gradient in the sense that there exists a Lyapunov function non-increasing along orbits, see \cite{Aragao-Costa-2011} for details.
\end{proof}

\begin{remark}\label{remark-reproductive-number-autonomous}
From Theorem \ref{th-SIR-autonomouos}, we can define now the \textbf{basic reproduction number}
\begin{equation}
    \mathcal{R}_0:=\frac{\gamma}{a+b+c},\label{R0}
\end{equation}
which represents the average number of susceptible individuals that an infected individual can infect during the infection period.
Thus, in the autonomous case, it is possible to describe the dynamics via $\mathcal{R}_0$ being less or greater than $1$. More precisely, the disease is eradicated as long as $\mathcal{R}_0<1$ and it becomes endemic (i.e., the number of infected individuals stays positive) whether $\mathcal{R}_0>1$.
We also note that \eqref{R0} can be also obtained from the well-known next-generation matrix method (see \cite{calculoR0,calculoR0old}).
\end{remark}

Since $T(\cdot)$ is gradient, it is stable under small nonautonomous and random perturbation, see \cite{Aragao-Costa-2013} and \cite{Caraballo-Carvalho-Oliveira-Sousa-NRA}, respectively. In the subsequent sections of this paper, we explore the impact of nonautonomous and random perturbations on the model, commencing our analysis with small disturbances.

\section{Small non-autonomous/random perturbations of a SIR model with vital dynamics and reinfection}
\label{SIR-small-nonautrandom-perturbations}

In this section, we consider small non-autonomous and random perturbations of the autonomous system \eqref{Snew}-\eqref{Rnew}.
First, we briefly recall the notions of random flow, non-autonomous random dynamical systems and attractors. In the sequel, $(X,d)$ denotes a complete separable metric space.
	
Let $(\Omega,\mathcal{F},\mathbb{P})$ be a probability space. 
We say that a family of maps $\{\theta_t:\Omega\rightarrow\Omega; \, t\in\mathbb{R}\}$ is a
\textbf{random flow} if 
	$\theta_0=Id_\Omega$;
		$\theta_{t+s}=\theta_t\circ \theta_s$, for all $t,s\in\mathbb{R}$;
		$\theta_t:\Omega \rightarrow \Omega$ is measurable and $\mathbb{P}\theta_t^{-1}=\mathbb{P}$ for all $t\in\mathbb{R}$.

	Define $\Theta_t(\tau,\omega):=(t+\tau,\theta_t\omega)$ for each $(\tau,\omega)\in\mathbb{R}\times\Omega$, and $t\in\mathbb{R}$.
	We say that a family of maps $\{\psi(t,\tau,\omega):X\to X; (t,\tau,\omega)\in\mathbb{R}^+\times\mathbb{R}\times\Omega\}$ 
	is a \textbf{nonautonomous random dynamical system (NRDS)} if it is a {\bf co-cycle} driven by
	$\Theta$, i.e., if
	\begin{enumerate}
		\item the mapping
		$\mathbb{R}^+\times\Omega\times X\ni (t,\omega,x)\mapsto \psi(t,\tau,\omega)x\in X$
		is measurable
		for each fixed $\tau\in\mathbb{R}$;
		\item
		$\psi(0,\tau,\omega)=Id_X$,
		for each $(\tau,\omega)\in\mathbb{R}\times\Omega$;
		\item 
		$\psi(t+s,\tau,\omega)=\psi(t,\Theta_s(\tau,\omega))\circ\psi(s,\tau,\omega)$,
		for every $t,s\geq 0$ in 
		$\mathbb{R}$, and $(\tau,\omega)\in\mathbb{R}\times \Omega$;
		\item $\psi(t,\tau,\omega):X\to X$ is a continuous map for each  $(t,\tau,\omega)\in\mathbb{R}^+\times\mathbb{R}\times\Omega$.
	\end{enumerate}
	We usually denote by $(\psi,\Theta)_{(X,\mathbb{R}\times\Omega)}$, or $(\psi,\Theta)$, the co-cycle $\psi$ driven by $\Theta$, and $\omega_\tau:=(\tau,\omega)\in \mathbb{R}\times \Omega$.
 
    \par Finally a {\bf random dynamical system (RDS)} is a co-cycle $\{\psi(t,\omega):X\to X; (t,\omega)\in\mathbb{R}^+\times\Omega\}$ driven by a random flow $\theta$, which we denote by $(\psi,\Theta)$.
    We say that $\eta:\mathbb{R}\times \Omega\to \mathbb{R}$ is a {\bf$\Theta$-invariant} function if $\eta(\Theta_t(\tau,\omega)=\eta(\Theta_t(\tau,\omega))$ for all $(\tau,\omega)\in \mathbb{R}\times \Omega$.

\par Now, we introduce the definitions of attractors and unstable sets for nonautonomous random dynamical systems. 
	Let $\hat{\mathcal{A}}=\{\mathcal{A}(\omega_\tau): \omega_{\tau}\in \mathbb{R}\times \Omega\}$ be a family of nonempty subsets of $X$. 
 
	We say that $\hat{\mathcal{A}}$ is a \textbf{nonautonomous random attractor} for $(\psi,\Theta)$
	if the following conditions are fulfilled:
	\begin{enumerate}
		
		\item $\mathcal{A}(\omega_{\tau})$ is compact, for every $\omega_\tau\in \mathbb{R}\times \Omega$;
		\item the set-valued mapping $\omega\mapsto \mathcal{A}(\tau,\omega)$ is measurable, for each $\tau \in \mathbb{R}$;
		\item $\hat{\mathcal{A}}$ is invariant, i.e., 
		$\psi(t,\omega_{\tau})\mathcal{A}(\omega_{\tau})=
		\mathcal{A}(\Theta_t\omega_{\tau})$ for every $t\geq 0$ and $\omega_{\tau}\in \mathbb{R}\times \Omega$;
		\item $\hat{\mathcal{A}}$ pullback attracts every bounded subset of $X$, i.e.,
		for every bounded subset $B$ of $X$ and $\omega_{\tau}\in \mathbb{R}\times\Omega$,
		\begin{equation*}
			\lim_{t\to +\infty}dist(\psi(t,\Theta_{-t}\omega_\tau)B, \mathcal{A}(\omega_{\tau}))=0,
		\end{equation*}
		where $\displaystyle dist(A,B)=\sup_{a\in A}\inf_{b\in B} d(a,b)$ is the Hausdorff semi-distance;
		\item $\hat{\mathcal{A}}$ is the minimal closed family that 
		pullback attracts bounded subsets of $X$. 
	\end{enumerate}
	

Subsequently, we employ the concept of unstable sets of global solutions for the nonautonomous random system introduced in \cite{Caraballo-Carvalho-Oliveira-Sousa-NRA}.

\par We say that $\xi:\mathbb{R}\times \Omega\to X$ is a \textbf{global solution} for the NRDS $(\psi,\Theta)$ if $\psi(t,\omega_\tau)\xi(\omega_\tau)=\xi(\Theta_t\omega_\tau)$.
	The \textbf{unstable set} of $\xi^*$ is the family
	\begin{equation*}
		W^u(\xi^*)=\{W^u(\xi^*;\omega_\tau): \omega_{\tau}\in \mathbb{R}\times \Omega\},
	\end{equation*}
	where, for each $\omega_{\tau}$, 
	$W^u(\xi^*;\omega_{\tau})$ is the unstable set of 
	the global solution $t\mapsto \xi^*(\Theta_t\omega_{\tau})$ of the evolution process
	$\Psi_{\omega_{\tau}}=\{\psi(t-s,\Theta_{s}\omega_{\tau}): t\geq s\}$, i.e.,
	\begin{eqnarray*}
		W^{u}(\xi^*;\omega_\tau)=\bigg\{(t,z)\in \mathbb{R}\times X: 	\hbox{ there is a global solution } \zeta \hbox{ of } 
		\Psi_{\omega_\tau} \hbox{ such that }\\
		\zeta(t)=z, \ 	\hbox{ and }
		\lim_{s\to -\infty} \|\zeta(s)-\xi^*(\Theta_s\omega_\tau)\|_X=0
		\bigg\}.
	\end{eqnarray*}
\par The authors in \cite{Caraballo-Carvalho-Oliveira-Sousa-NRA} provide conditions under which the gradient dynamics is preserved under nonautonomous random perturbations, such that the perturbed attractors are equal to the union of the unstable sets of the random hyperbolic solutions. Moreover, the authors in \cite{Caraballo-Carvalho-Langa-OliveiraSousa-2021} study the notion of random hyperbolic solutions and exponential dichotomies for NRDS.

\par We consider small nonautonomous random perturbations of the autonomous case \eqref{Snew}-\eqref{Rnew} to apply the abstract results of \cite{Caraballo-Carvalho-Oliveira-Sousa-NRA}. Our goal is to prove that, even with these types of perturbations, it is possible to characterize the attractors and interpret if the disease is eradicated or not.

\par Let $(\Omega, \mathcal{F},\mathbb{P})$ be a probability space and
$\theta$ be a random flow defined in $\Omega$. 
Consider a family of perturbations $F_\eta(t,\theta_t\omega,S,I,R)\in \mathbb{R}^3$, with $\eta\in(0,1]$, of \eqref{Snew}-\eqref{Rnew}, that depends on the time parameter $t$ and the random parameter $\omega$. We obtain the following results for this type of perturbation directly from  \cite[Theorem 5.1 and Theorem 6.1]{Caraballo-Carvalho-Oliveira-Sousa-NRA}.

\begin{corollary}\label{th-topological-structure-stability} 
Let $\gamma$ be such that $\gamma>a+b+c$ and 
that $F_\eta$ is such that
\begin{equation}\label{eq-condition-F_eta}
    \sup_{t\in \mathbb{R}}\|F_\eta(t+\tau,\theta_t\omega,\cdot)\|_{C^1(\mathbb{R}^3)}\to 0, \ \forall \, (\tau,\omega)\in \mathbb{R}\times \Omega.
\end{equation}
We also assume that $F_\eta$ is such that
the perturbed problem is globally well-posed and associated with a co-cycle $(\psi_\eta,\Theta)$. Furthermore, assume that
for each $\eta\in [0,1]$, the co-cycle $(\psi_\eta,\Theta)$ has a nonautonomous random attractor $\{\mathcal{A}_\eta(\omega_\tau): \omega_\tau\in \mathbb{R}\times \Omega\}$,
	   \begin{equation*}\label{eq-continuity-at-epsilon-0-thcontinuityattractors}
			K(\omega_{\tau}):=\overline{\bigcup_{t\in \mathbb{R}} \bigcup_{\eta\in [0,1]} \mathcal{A}_\eta(\Theta_t\omega_\tau)} \hbox{ is compact,} \ \ \forall \, \omega_{\tau}\in \mathbb{R}\times \Omega, \hbox{ and }
		\end{equation*}
		\begin{equation}\label{th-continuity-at-epsilon-0-thcontinuityattractors-2}
		\bigcup_{\eta\in [0,1]} \bigcup_{t\geq 0} \psi_{\eta}(t,\omega_{\tau})	K(\omega_{\tau}) \hbox{ is bounded,} \ \ \forall \, \omega_{\tau}\in \mathbb{R}\times \Omega.
		\end{equation}
Then the nonautonomous random dynamical systems associated with
the perturbed equations have continuity and topological structural stability. 
More precisely, given $\epsilon_0>0$ suitable small, there exists a $\Theta$-invariant function $\eta_{0}:\mathbb{R}\times \Omega\to (0,1]$ such that, for each $\omega_{\tau}$ fixed, the following hold:
	\begin{enumerate}
		\item For any $\eta \in (0,\eta_0(\omega_{\tau})]$ and $j\in \{0,1\}$, there exists a hyperbolic solution $\xi_{j,\eta}^*$ of $\{\psi_\eta(t-s,\Theta_s\omega_{\tau}): t\geq s\}$ 
		with
		\begin{equation}
		\sup_{j} \sup_{t\in \mathbb{R}} \|\xi_{j,\eta}^*(\Theta_{t}\omega_\tau)-E_j^*\|_X<\epsilon_0,
		\end{equation}
		where 
		the linearized associated evolution process  
		admits an exponential dichotomy with family of projections $\{\Pi_{j,\eta}^u(s): s\in \mathbb{R}\}$.

		\item For each $\omega_{\tau}$, the family of pullback attractors
		$\{\mathcal{A}_\eta(\Theta_t\omega_{\tau}): t\in \mathbb{R}\}_{\eta\in [0,\eta_{0}(\omega_{\tau})]}$ is continuous at $\eta=0$.
	In particular, we have continuity of nonautonomous random attractors in the following sense: 
	given $\epsilon>0$, there exists a $\Theta$-invariant function
	$\eta_{\epsilon}\leq \eta_0$ such that,  
	for every $\Theta$-invariant function $\bar{\eta}$, with
	$\bar{\eta}\leq \eta_{\epsilon}$, we have
	\begin{equation}
	\sup_{t\in \mathbb{R}}d_H(\mathcal{A}_{\bar{\eta}}(\Theta_t\omega_{\tau}),\mathcal{A}_0 )<\epsilon, \ \ 
	\forall \, \omega_{\tau}\in \mathbb{R}\times \Omega,
	\end{equation}
	where $\{\mathcal{A}_{\bar{\eta}}(\omega_{\tau}): \omega_{\tau}\in \mathbb{R}\times \Omega\}$ is the nonautonomous random attractor of 
	$(\psi_{\bar{\eta}},\Theta)$ and 
	$d_H(A,B)=\max\{dist_H(A,B),dist_H(B,A)\}$, for $A,B\subset X$.
	\item Then, 
	there exists a $\Theta$-invariant function $\eta_{1}:\mathbb{R}\times \Omega\to (0,1)$ 
	such that for each $\omega_{\tau}$ fixed
the evolution process
$\{\psi_\eta(t-s,\Theta_{s}\omega_{\tau}): t\geq s\}$ is
dynamically gradient with respect to $\{\xi_{0,\eta}^*,\xi_{1,\eta}^*\}$, $\forall \eta \leq \eta
_1(\omega_{\tau}) $. 
Consequently, 
\begin{equation}
\mathcal{A}_\eta(\Theta_t\omega_\tau)= \{\xi_{1,\eta}(\Theta_t\omega_\tau)\}\cup W^u_\eta(\xi_{0,\eta}^*;\omega_{\tau})(t),  \forall\, \eta\in [0,\eta_1(\omega_{\tau})].
\end{equation}
\end{enumerate}
\end{corollary}
The following result summarizes the case where the disease-free equilibrium is the only hyperbolic equilibrium of the autonomous systems. 

\begin{corollary}\label{th-topological-structure-stability-2} 
    Suppose that $\gamma<a+b+c$ and that $F_\eta$ satisfies the conditions of Corollary \ref{th-topological-structure-stability}. Given
     $\epsilon_0>0$ suitable small,
    then 
	there exists a $\Theta$-invariant function $\eta_{1}:\mathbb{R}\times \Omega\to (0,1)$ 
	such that for each $\omega_{\tau}$ fixed and $\eta\leq \eta_1(\omega_\tau)$ there exists a hyperbolic solution $t
\mapsto \xi_{0,\eta}^*(\Theta_t\omega_\tau)$ for
the evolution process 
$\{\psi_\eta(t-s,\Theta_{s}\omega_{\tau}): t\geq s\}$	with
		\begin{equation}
		\sup_{t\in \mathbb{R}} \|\xi_{0,\eta}^*(\Theta_{t}\omega_\tau)-E_0^*\|_X<\epsilon_0.
		\end{equation}
Moreover, $\{\psi_\eta(t-s,\Theta_{s}\omega_{\tau}): t\geq s\}$
dynamically gradient with respect to $\{\xi_{0,\eta}^*(\Theta_{(\cdot)}\omega_\tau\}$, $\forall \eta \leq \eta
_1(\omega_{\tau}) $. 
Consequently, 
$\mathcal{A}_\eta(\omega_\tau)= \{\xi_{0,\eta}(\omega_\tau)\},  \   \omega_\tau\in \mathbb{R}\times \Omega.$ 
\end{corollary}

As an application, one can perturb the transmission coefficient, considering 
$\gamma_\eta(\theta_t\omega_\tau)=\gamma+\eta \phi(z^*(\theta_t))$ instead of $\gamma$ in \eqref{Snew} and \eqref{Inew}, where $\phi$ is a bounded function. Then the associated perturbation $F_\eta$ also satisfies the condition of Corollary \ref{th-topological-structure-stability} and it is straightforward to verify Condition \eqref{eq-continuity-at-epsilon-0-thcontinuityattractors} for the perturbed problems. 

Another possible application is to consider
the following stochastic environmental  perturbation 
\begin{align*}
dS&=(q-aS+bI)dt-\gamma\frac{SI}{S+I+R}dt+\eta \kappa(t)S \circ dW_t,\\
dI&=-(a+b+c)Idt+\gamma\frac{SI}{S+I+R}dt +\eta \kappa(t)I \circ dW_t,\\
dR&=cIdt-aR dt +\eta \kappa(t)R \circ dW_t,\end{align*}
where $\{W_t\}_{t\in \mathbb{R}}$  is the Wierner process, $\eta>0$, $\kappa:\mathbb{R}\to \mathbb{R}$ is a continuously differentiable function with compact support in $\mathbb{R}$ or behaves as $1/(1+|t|)$ for $|t|$ suitable large, 
                   see \cite[page 23]{Caraballo-Carvalho-Oliveira-Sousa-NRA} for details. 

Then, in both classes of examples above one can obtain continuity of attractors (upper and lower semicontinuity of the attractors) and topological structural stability for the above system applying Corollary \ref{th-topological-structure-stability} and Corollary
\ref{th-topological-structure-stability-2}.

\begin{remark}\label{remark_3.1}
    The results of this section can be understood as follows. If $\mathcal{R}_0>1$, the disease becomes endemic. Conversely, when $\mathcal{R}_0<1$, the number of infected individuals approaches zero as time tends to infinity. This behavior aligns with the reproductive number $\mathcal{R}_0$ in the autonomous scenario, as discussed in Remark \ref{remark-reproductive-number-autonomous}.
\end{remark}

In the following sections, we consider random perturbations (not necessarily small) of some parameters in system \eqref{Snew}-\eqref{Rnew}, motivated by real phenomena (see, for instance \cite{Albani-Zubelli}). In each case, we provide some information about the internal structure of the corresponding pullback random attractor and the biological interpretation. 

\section{SIR model with vital dynamics, reinfection and random fluctuations on the diffusion coefficient}\label{SIRgamma}

We study the SIR model with vital dynamics, reinfection, and random fluctuations on the transmission coefficient, given by
\begin{align}
S^\prime&=q-aS+bI-(\gamma+\Phi(\theta_t\omega))\frac{SI}{S+I+R}\label{odeS}\\
I^\prime&=-(a+b+c)I+(\gamma+\Phi(\theta_t\omega))\frac{SI}{S+I+R}\label{odeI}\\
R^\prime&=cI-a,\label{odeR}
\end{align}
where $S$, $I$ and $R$ describe again the number of susceptible, infected and recovered individuals, respectively, and the parameters have already been introduced in Section \ref{introduction}.

\begin{theorem}\label{th-SIR-perturbing-gamma}

For any $\omega\in\Omega$ and $u_0
\in\R^3_+$, system \eqref{odeS}-\eqref{odeR} has a unique solution $u(\cdot,\omega,u_0)
\in  {C}^1([0,+\infty);\R^3_+)$, which generates a RDS $\varphi:[0,+\infty)\times\Omega\times\R^3_+\rightarrow\R^3_+$ defined as $\varphi(t,\omega)u_0=u(t;\omega,u_0)$ for all $t\in[0,+\infty)$, $u_0\in\R^3_+$ and $\omega\in\Omega$. In addition, 
%
there exists a unique pullback random attractor $\mathcal{A}=\{\mathcal{A}(\omega)\}_{\omega\in\Omega}$ for the RDS $(\varphi,\theta)$, whose component subsets satisfy
$\mathcal{A}(\omega)\subset B_0$, where
\begin{equation}
B_0:=\left\{(S,I,R)\in\R^3_+\,:\,S+I+R=\f{q}{a}\right\}.\label{B0qconstante}
\end{equation}
\end{theorem}

\begin{proof}
First, we fix $\omega\in \Omega$. The local well-posedness of \eqref{odeS}-\eqref{odeR} is standard. The first step is to prove that $\R^3_+$ is positively invariant for the solutions of system \eqref{odeS}-\eqref{odeR}, i.e., every local solution with an initial condition in $\R^3_+$ remains in $\R^3_+$ for every time. 

By continuity, if a component ($S$, $I$ or $R$) of any solution $u$ of system \eqref{odeS}-\eqref{odeR} with initial condition in $\R^3_+$ reaches a negative value, then there exists some time $t^*\in\R$ at which that component ($S$, $I$ or $R$) takes the value zero. 
Assume, for instance, that $S(t^*)=0$, $I(t^*)\geq 0$ and $R(t^*)\geq 0$ for some $t^*\in\R$. From \eqref{odeS}, we have that $S'(t^*)=q+bI(t^*)>0$, whence we deduce that the plane $S=0$ cannot be crossed by any solution of system \eqref{odeS}-\eqref{odeR} with an initial condition in $\R^3_+$. By a similar reasoning, $R$ cannot take negative values.
Finally, since every solution starting in the plane $I=0$ remains there for every future time, then $I$ cannot take negative values either. Thus, $\R^3_+$ is positively invariant.

Let us define $N=S+I+R$, which satisfies $N'=q-aN$. Hence
\begin{equation}
    N(t,N_0)=e^{-at}N_0+\f{q}{a}\left[1-e^{-at}\right],\label{solN}
\end{equation}
\noindent for every $t\geq 0$, $\omega\in\Omega$ and $N_0=S_0+I_0+R_0\geq 0$.

From \eqref{solN}, it is clear that $N$ does not blow up at any finite time. Thus, since $S$, $I$ and $R$ remain non-negative for every $t\in\R$, $S$, $I$ and $R$ do not blow up at any finite time either, whence every unique local solution of system \eqref{odeS}-\eqref{odeR} is globally defined in time.

In view of the previous arguments, it is straightforward to verify that the solution mapping $u(\cdot,\omega,u_0)$ of system \eqref{odeS}-\eqref{odeR} generates a RDS $\varphi:\mathbb{R}_+\times\Omega\times\R^3_+\rightarrow\R^3_+$ defined as $\varphi(t,\omega)u_0=u(t;\omega,u_0)$ for all $t\in\R$, $u_0\in\R^3_+$ and $\omega\in\Omega$.

Then, from \eqref{solN},
\begin{equation}
    \lim_{t\rightarrow+\infty}N(t,N_0)=\f{q}{a},\label{limNq/a}
\end{equation}
for $N_0=S_0+I_0+R_0\geq 0$ uniformly in $\Omega$.
Actually, for any $u_0\in \mathbb{R}^3_+$ we have
\begin{equation}\label{eq-uniform-attraction-B_0}
    \lim_{t\to +\infty} \sup_{\omega\in \Omega} d(\varphi(t,\omega)u_0,B_0)=0,
\end{equation}
in particular, $d(\varphi(t,\theta_{-t}\omega)u_0,B_0)\to 0$ as $t\to +\infty$, so 
$B_0$ attracts pullback bounded subsets of $\mathbb{R}^3_+$.
In this case, it is well-known that there exists a unique pullback random attractor $\mathcal{A}=\{\mathcal{A}(\omega)\}_{\omega\in\Omega}$ for the RDS $\varphi$, with $\mathcal{A}(\omega)\subset B_0$.
\end{proof}

\begin{remark}
Note that, from the proof of the theorem, the absorbing set is deterministic, i.e., it does not depend on the realization of the noise. Note also that \eqref{eq-uniform-attraction-B_0} is a forward and uniform (with respect to $ \omega$) attraction.
\end{remark}

The perturbation on the transmission parameter $\gamma$ does not affect the existence of attractors because it does not appear in the equation of $N$, however, as we will see in Theorem \ref{th-existence-xi*-hyperbolicity}, it has an effect on the structure of the attractor.

\begin{theorem}\label{th-q-cte-gamma-pertturbed}
Let $\gamma,a,b,c>0$ be such that $\gamma<a+b+c$
and assume that
\begin{equation}\label{eq-ergodic-condition-for-Phi}
	\lim_{t\to \pm\infty} \frac{1}{t}\int_{0}^{t}\Phi(\theta_r\omega)\, dr=0, \ \hbox{ for } \omega\in \Omega.
	\end{equation}
Then the RDS associated to \eqref{odeS}-\eqref{odeR} has a singleton random attractor $\mathcal{A}$ given by $\mathcal{A}(\omega)=E_0^*=(q/a,0,0)$ for every $\omega\in\Omega$, which means that the disease is eradicated.
\end{theorem}

\begin{proof}
Let $u_0=(S_0,I_0,R_0)$ be such that $S_0,I_0,R_0>0$.
From \eqref{odeI}, since $S/N\leq 1$, 
\begin{equation}
\f{dI}{dt}\leq\left[-(a+b+c)+\gamma+\Phi(\theta_t\omega))\right]I.\label{Isol0}
\end{equation}
Then $
I(t,\omega,u_0)\leq I_0e^{-(a+b+c-\gamma)t-\int_{0}^t\Phi(\theta_r\omega)dr}$, for every $t\geq 0$, $\omega\in\Omega$ and $I_0>0$.

Thus, since $\gamma<a+b+c$ and \eqref{eq-ergodic-condition-for-Phi} fulfills,  replacing $\omega$ by $\theta_{-t}\omega$ and taking limit
\begin{equation}\lim_{t\rightarrow+\infty}I(t,\theta_{-t}\omega,u_0)=0, 
 \ \forall \, \omega\in\Omega.\label{limitI0}
\end{equation}
Let $\omega\in \Omega$ and $\eta>0$ arbitrary. From \eqref{limitI0}, there exists $T_0(\omega,\eta)>0$ such that $I(t,\theta_{-t}\omega,u_0)\leq\eta$ for all $t\geq T_0(\omega,\eta)$. Then, from \eqref{odeR}, we have
  $  R^\prime\leq c\eta-aR,$
for every $t\geq T_0(\omega,\eta)$. Thus
\begin{equation}
    R(t,\theta_{-t}\omega,u_0)\leq R(t_0,\theta_{-t_0}\omega,u_0) e^{-a(t-t_0)}+\f{c\eta}{a}\left[1-e^{-a(t-t_0)}\right], 
\end{equation}
for $t\geq t_0\geq  T_0(\omega,\eta)$.
Then there exists $T_1(\omega,\eta)\geq \max\{t_0,T_0(\omega,\eta)\}$ such that 
\begin{equation}
    R(t,\theta_{-t}\omega,u_0)\leq \eta(1+\f{c}{a}), \ \hbox{ for } t\geq T_1(\omega,\eta).
\end{equation}
\noindent Hence,
\begin{equation}
 \lim_{t\rightarrow+\infty} R(t,\theta_{-t}\omega,u_0)=0, \ \forall \, \omega\in\Omega.
\end{equation}

Then, since both $I(t,\theta_{-t}\omega,u_0)$ and $R(t,\theta_{-t}\omega,u_0)$ converge to $0$ as $t\to +\infty$, from \eqref{limNq/a} we obtain
\begin{equation*}
    \lim_{t\rightarrow+\infty}S(t,\theta_{-t}\omega,u_0)=\f{q}{a}, \ \forall\, \omega\in\Omega.
\end{equation*}

Finally, since from Theorem \ref{th-SIR-perturbing-gamma} $(\varphi,\theta)$ has a pullback attractor $\mathcal{A}=\{\mathcal{A}(\omega)\}_{\omega\in\Omega}$ and $E_0^*$ 
pullback attracts every bounded subset of $\mathbb{R}^3_+$, we see that $\mathcal{A}(\omega)=\{E_0^*\}$ for every $\omega\in\Omega$.
\end{proof}


\par Now, we provide conditions under which the disease becomes endemic.

\begin{theorem}\label{th-condition-below-boundness-gamma}
Let $\gamma_0>0$ be such that 
$\gamma_0\leq \gamma+\Phi(\omega)$ for all $\omega\in \Omega$.
If $\gamma_0>a+b+c$,
then $I$ remains positive forward and pullback asymptotically, i.e., for all $u_0\in \mathbb{R}^3_+$ with $I_0>0$,
there exists $\epsilon_0>0$ such that
\begin{equation}
        \lim_{t\to +\infty}I(t,\omega,u_0)>\epsilon_0 \hbox{ and} 
    \lim_{t\to +\infty}I(t,\theta_{-t}\omega,u_0)>\epsilon_0.
\end{equation}
\end{theorem}

\begin{proof}

We are going to prove the first limit, the other is analogous.
From \eqref{limNq/a}, for any $\eta>0$ there exists $t_0=t_0(\eta)>0$ such that
$$\f{q}{a}-\eta\leq N(t,N_0)\leq\f{q}{a}+\eta,$$
\noindent for all $t\geq t_0$ and $N_0=S_0+I_0+R_0>0$.

Consider now $\varepsilon\in(0,I(t_0,\omega,u_0))$ and assume that
$I(t,\omega,u_0)<\varepsilon$ for all $t\geq t_0$.
Then, from \eqref{odeR}, we deduce
$R(t,\omega,u_0)\leq c\varepsilon/{a}, \ t\geq t_0$.
Since $N=S+I+R$, 
\begin{equation}
S(t,\omega,u_0)\geq N(t,N_0)-\varepsilon-\f{c\varepsilon}{a}, \ t\geq t_0.
\end{equation}

Thus, from \eqref{odeI}, for $t\geq t_0,$
\begin{align}
\nonumber
I^\prime(t,\omega,u_0)&\geq\left[-(a+b+c)+(\gamma+\Phi(\theta_t\omega))\f{N(t,N_0)-\varepsilon-\f{c\varepsilon}{a}}{N(t,N_0)}\right]I(t,\omega,u_0)\\
&\geq\left[-(a+b+c)+\gamma_0\f{\f{q}{a}-\eta-\varepsilon-\f{c\varepsilon}{a}}{q+\eta}\right]I(t,\omega,u_0).\label{I>0qconstante}
\end{align}
Now since $\gamma_0>a+b+c$, 
\begin{equation*}
-(a+b+c)+\gamma_0\f{\f{q}{a}-\eta-\varepsilon-\f{c\varepsilon}{a}}{q+\eta}>0
\end{equation*}
\noindent for $\eta>0$ and $\varepsilon>0$ small enough.

Hence, from \eqref{I>0qconstante}, $I^\prime(t_0,\omega,u_0)>0$, which completes the proof.
\end{proof}


\begin{remark}\label{Remark_reproductive_number_fails}
From Theorem \ref{th-q-cte-gamma-pertturbed} we can define $\mathcal{R}_2=\gamma/(a+b+c)$ such that if $\mathcal{R}_2<1$ the disease will be eradicated. And given  $\mathcal{R}_1=\gamma_0/(a+b+c)$, from Theorem \ref{th-condition-below-boundness-gamma} the disease becomes endemic if $\mathcal{R}_1>1$. We do not know if Theorem \ref{th-condition-below-boundness-gamma} holds true for $\gamma>a+b+c$ such that $\mathcal{R}_1=\mathcal{R}_2$, but it will be interesting to have an analogous reproductive number of the autonomous case (see Remark \ref{remark-reproductive-number-autonomous}) to the random SIR model. 
\end{remark}

Using the language of \cite{Kloeden-Kozyakin} and \cite{Kloeden-Potzsche}, Theorem \ref{th-condition-below-boundness-gamma} says that the disease-free equilibrium 
$S^*=(q/a,0,0)$ is unstable. In the following section, we obtain that this solution is hyperbolic with an one-dimensional unstable manifold (see Theorem \ref{theorem_cases_q_no_constant} with constant $q$) also in the case where $q$ is not a constant, similar to the autonomous case (Lemma \ref{lemma-E_0-E_1-hyperbolics}).

Note that the random variable $\Phi$ could be considered nonautonomous and random, $\Phi=\Phi(t,\omega)$, and we could prove similar results for the NRDS associated. Taking this into consideration, next we explore a recruitment rate, denoted by $q$, that is both nonautonomous and random.

\section{SIR model with vital dynamics, reinfection and random fluctuations on the
diffusion coefficient and the recruitment rate}\label{SIRqgamma}



For this section, we assume that $q:\mathbb{R}\times\Omega\to (0,+\infty)$ is such that
\begin{enumerate}
\item $q(\tau,\cdot)$ is a random variable for every $\tau$;
\item there are $\Theta$-invariant functions $q_0,q_1:\mathbb{R}\times \Omega\to (0,+\infty)$ such that $q_0(\omega_\tau)\leq q(\omega_\tau)\leq q_1(\omega_\tau)$;
\item $t\mapsto q(t,\theta_t\omega)$ is continuous for each $\omega \in \Omega$.
\end{enumerate}

Now, we analyze system \eqref{S_random_gamma}-\eqref{R_random_gamma} with the above recruitment rate $q$. By a translation in time of \eqref{S_random_gamma}-\eqref{R_random_gamma}, see \cite[Sectiton 3]{Caraballo-Carvalho-Langa-OliveiraSousa-2021}, we study 
\begin{align}
S^\prime&=q(t+\tau,\theta_t\omega)-aS+bI-(\gamma+\Phi(\theta_t\omega))\frac{SI}{S+I+R},\label{eS1}\\
I^\prime&=-(a+b+c)I+(\gamma+\Phi(\theta_t\omega))\frac{SI}{S+I+R},\label{eI1}\\
R^\prime&=cI-aR,\label{eR1}
\end{align}
for $t\geq 0$, $\tau \in \mathbb{R}$ and $\omega\in \Omega$. 
We prove the existence of attractors and random hyperbolic solutions for the NRDS associated with \eqref{eS1}-\eqref{eR1}. In addition, we provide conditions such that the disease is eradicated or not.

\begin{theorem}
For any $\omega\in\Omega$, $\tau\in \mathbb{R}$ and $u_0\in\R^3_+$, system \eqref{eS1}-\eqref{eR1} has a unique global solution $u(\cdot,\tau,\omega,u_0)\in C^1([0,+\infty));\mathbb{R}^3_+)$, 
which generates a NRDS $\varphi:[0,+\infty)\times \mathbb{R}\times\Omega\times\R^3_+\rightarrow\R^3_+$ defined as $\varphi(t,\omega_\tau)u_0=u(t,\tau,\omega,u_0)$ for all $t\geq 0$.
Moreover, there exists a nonautonomous random attractor $\mathcal{A}=\{\mathcal{A}(\omega_\tau)\}_{\omega_\tau\in\mathbb{R}\times\Omega}$ for the NRDS $(\varphi,\Theta)$, such that
$$\mathcal{A}(\omega_\tau)\subset B_0(\omega_\tau)=\left\{(S,I,R)\in\R^3_+\,:\,\f{q_0(\omega_\tau)}{a}\leq S+I+R\leq \f{q_1(\omega_\tau)}{a}\right\},  \ \omega_\tau\in\mathbb{R}\times\Omega.$$
\end{theorem}

\begin{proof}
The proof of local well-posedness and that $\mathbb{R}^3_+$ is invariant is similar to the proof of Theorem \ref{th-SIR-perturbing-gamma}.
Now we prove that the solutions are globally defined. Note that 
 $N=S+I+R$ satisfies
$    N^\prime=q(\Theta_t\omega_\tau)-aN.$
Thus 
\begin{equation}
    N(t,\omega_\tau;N_0)=e^{-at}N_0+\int_0^t e^{-a(t-s)} q(\Theta_s\omega_\tau)ds, \ t>0.
\end{equation}
Then it follows that the solutions are globally defined. From the soundness of $q$, we have
\begin{equation}
    e^{-at}N_0+\f{q_0(\omega_\tau)}{a}\left[1-e^{-at}\right]\leq N(t,\omega_\tau,N_0)
    \leq    e^{-at}N_0+\f{q_1(\omega_\tau)}{a}\left[1-e^{-at}\right],
\end{equation}
for every $t\geq 0$ and $\omega_\tau$. Therefore $\{B_0(\omega_\tau):\omega_\tau\in \mathbb{R}\times \Omega\}$ pullback (and forward) attracts bounded subsets of $\mathbb{R}^3_+$ and this concludes the proof.
\end{proof}

The subsequent outcome arises inspired in Theorem \ref{th-q-cte-gamma-pertturbed} and Theorem \ref{th-condition-below-boundness-gamma}, and the some parts of the will be omitted.

\begin{theorem}\label{theorem_cases_q_no_constant}
The following statements hold true:
\begin{enumerate}
\item If $\gamma<a+b+c$ and \eqref{eq-ergodic-condition-for-Phi} is satisfied, then  $\displaystyle{\lim_{t\rightarrow+\infty}I(t,\Theta_{-t}\omega_\tau,u_0)=0}$. Actually, there here exist a global solution $\xi^*(\omega_\tau):=(S^*(\omega_\tau),0,0)$ such that 
$\mathcal{A}(\omega_\tau)=\{\xi^*(\omega_\tau)\}$ for $ \omega_\tau\in \mathbb{R}\times \Omega$.
\item Assume that $q_0$ and $q_1$ are constants. If $\gamma_0 q_0/q_1>a+b+c$, then 
for all $u_0=(S_0,I_0,R_0)\in \mathbb{R}^3_+$ with $I_0>0$
there exists $\epsilon_0>0$ such that
\begin{equation}
        \lim_{t\to +\infty}I(t,\omega_\tau,u_0)>\epsilon_0 \hbox{ and} 
    \lim_{t\to +\infty}I(t,\Theta_{-t}\omega_\tau,u_0)>\epsilon_0.
\end{equation}
\end{enumerate}

\end{theorem}

\begin{proof}
\par We only prove the first item, the remainder of the proof follows is similar to the mentioned results and will be omitted.
The proof of the first statement is similar to the proof of Theorem \ref{th-q-cte-gamma-pertturbed} and is omitted.
For the second claim $I(t,\Theta_{-t}\omega_\tau)u_0\to 0$, we have $R(t,\Theta_{-t}\omega_\tau)u_0\to 0$. If $I=R=0$, the equation for $S$ is equivalent to 
\begin{equation}
    S(t,\omega_\tau,u_0)=e^{-at}S_0+\int_0^t e^{-a(t-s)}q(\Theta_s\omega_\tau) ds.
\end{equation}
Then, we define
\begin{equation}
        S^*(\omega_\tau):=\lim_{t\to +\infty}S(t,\Theta_{-t}\omega_\tau,u_0)
        = \int_0^{+\infty} e^{-as}q(\Theta_{-s}\omega_\tau) ds.
\end{equation}
Note that $\xi^*(\omega_\tau):=(S^*(\omega_\tau),0,0)$  
is a global solution for $(\varphi,\Theta)$
that pullback attracts bounded subsets of $\mathbb{R}^3_+$ and 
the nonautonomous random attractor is given by 
$\mathcal{A}(\omega_\tau)=\{\xi(\omega_\tau)\}$.
The remainder of the proof follows is similar to the proof of Theorem \ref{th-condition-below-boundness-gamma} and will be ommited.
\end{proof}
\begin{remark}
    The condition to ensure the extinction of the disease in \cite{Caraballo-Colucci-17} depends on the amplitude of the random perturbations, whereas our condition does not involve the noise, i.e., it is finer thanks to the way in which we model the noise.
\end{remark}

The global solution $\xi^*(\omega_\tau):=(S^*(\omega_\tau),0,0)$, which represents the deasease-free ``nonautonomous equilibrium''. In \cite{Kloeden-Kozyakin} they say that 
$\xi^*$ is a unstable ``non-autonomous equilibirum'' if $I$ remains positive pullback at infinity. The following result is our attempt to understand more about hyperbolicity of this special solution.

\begin{theorem}\label{th-existence-xi*-hyperbolicity}
Let $\xi^*$ be defined as by $\xi^*(\omega_\tau):=(S^*(\omega_\tau),0,0)$. Then $\xi^*$ is a global solution for $(\varphi,\Theta)$ satisfying the following:
\begin{enumerate}
\item if $\gamma<a+b+c$, and $\Phi$ satisfies \eqref{eq-ergodic-condition-for-Phi}, then $\xi^*$ is a random hyperbolic solution: the linearized problem around it admits a tempered exponential decay. In particular, if $\Phi=0$, then $\xi^*$ is asymptotically exponentially stable. 
\item if $\gamma>a+b+c$ and $\Phi=0$, then the global solution $\xi^*$ is a random hyperbolic solution with a one-dimensional local unstable manifold, i.e., the co-cycle associated with the linearization around $\xi^*$ 
admits a uniform exponential dichotomy with a one-dimensional unstable projection;
\end{enumerate}
\end{theorem}

\begin{proof}
The proof of the existence of $\xi^*$ is similar to the proof of Theorem \ref{theorem_cases_q_no_constant}, thus we prove that $\xi^*$ is hyperbolic. 

First, we assume that $\Phi=0$, then the linearization of \eqref{eS1}-\eqref{eR1} around $\xi^*$ provides us the following matrix
\begin{equation}\label{eq-A-case-phi=0}
A=\left(\begin{array}{ccc}
-a & b-\gamma & 0\\
0 & \gamma-a-b-c & 0\\
0 & c & -a
\end{array}\right),
\end{equation}
which is straightforward to see that admits exponential dichotomy by taking the projections into respective eigenspaces of $-a$ and $\lambda_1=\gamma-a-b-c$. 

Next, we prove item 1. Assume that $\Phi$ satisfies \eqref{eq-ergodic-condition-for-Phi}. The linearized problem around $\xi^*$ is given by 
$\dot{x}=Ax+B(\theta_t\omega)x$ and the co-cycle associated is
$
   \varphi(t,\omega)=e^{At+\int_0^tB(\theta_s\omega)ds},$
where $A$ is given in \eqref{eq-A-case-phi=0} and 
\begin{equation*}
B(\omega)=\Phi(\omega)\left(\begin{array}{ccc}
0 & 1 & 0\\
0 & 1 & 0\\
0 & 0 & 0
\end{array}\right), \ \omega\in \Omega.
\end{equation*}

Thus
 $   \|\varphi(t,\omega)\|\leq K(t,\omega)e^{-\lambda t},$ 
 where $K(t,\omega)=e^{\int_0^t\Phi(\theta_s\omega)ds}$ is tempered for each $\omega$ fixed, i.e., 
\begin{equation*}
    \frac{\ln K(t,\omega)}{t}\to 0, \hbox{ as } t\to +\infty.
\end{equation*}
Therefore, $\xi^*$ is a tempered hyperbolic solution and the proof is finished.
\end{proof}

\begin{remark}\label{Remark_5.1}
In our opinion, the case $\gamma>a+b+c$ and $\Phi$ not zero is also possible to study, even though it is an open question up to now. One alternative is to use \cite[Theorem 3.2]{Barreira-Valls-upper-tempered} or to follow the ideas of \cite{Battateli-Palmer-2015} for the case where $\Phi$ is bounded. 

We also note that, if $\mathcal{R}_0<1$ and $\Phi$ satisfies \eqref{eq-ergodic-condition-for-Phi}, the disease is eradicated. Thus it remains open the question: if $\mathcal{R}_0>1$ and $\Phi$ satisfying \eqref{eq-ergodic-condition-for-Phi} the disease became endemic?
\end{remark}

Some examples of the recruitment rate $q$:
\begin{enumerate}
    \item Define $q(t,\omega)=\nu e^{z(\omega)\kappa(t)}$, where
$\nu >0$, $z$ is the Ornstein-Uhlenbeck process, and 
$\kappa:\mathbb{R}\to \mathbb{R}$ is a continuously differentiable function with compact support in $\mathbb{R}$ or behavesr as $1/(1+|t|)$ for $|t|$ large enough, see \cite[page 23]{Caraballo-Carvalho-Oliveira-Sousa-NRA} for more details of this type of examples.
\item $q(t,\omega)=q(t)$ be 
any bounded time-dependent positive function. So in this case, the above results are reinterpreted with pullback attractors and hyperbolic solutions for the associated evolution process, see \cite{Carvalho-Langa-Robison-book}.
\end{enumerate}

\section{Numerical simulations}\label{sec:numerical_simulations}

In this section, we present several numerical simulations to illustrate the results proved throughout our paper. In the sequel, every figure shows a big panel on the right, depicting the phase plane of the corresponding system with an arrow pointing at the initial condition, and three little panels on the left, where we display the evolution on time of the number of susceptible, infected and recovered individuals. In addition, in each panel the blue dashed/dotted line represents the solution of the deterministic systems (i.e., with no noise) and the rest colored continuous lines are different realizations of the random models.

\subsection{Autonomous SIR model with vital dynamics and reinfection}

Our aim in this subsection is to display numerical simulations concerning the deterministic autonomous SIR model with reinfection \eqref{Snew}-\eqref{Rnew}. To this end, we set $a=1.5$, $b=0.5$, $c=0.7$, $q=5$ and we consider the initial condition $(S_0,I_0,R_0)=(25,1,0)$.

\begin{figure}[H]
    \begin{center}
    \includegraphics[width=\textwidth]{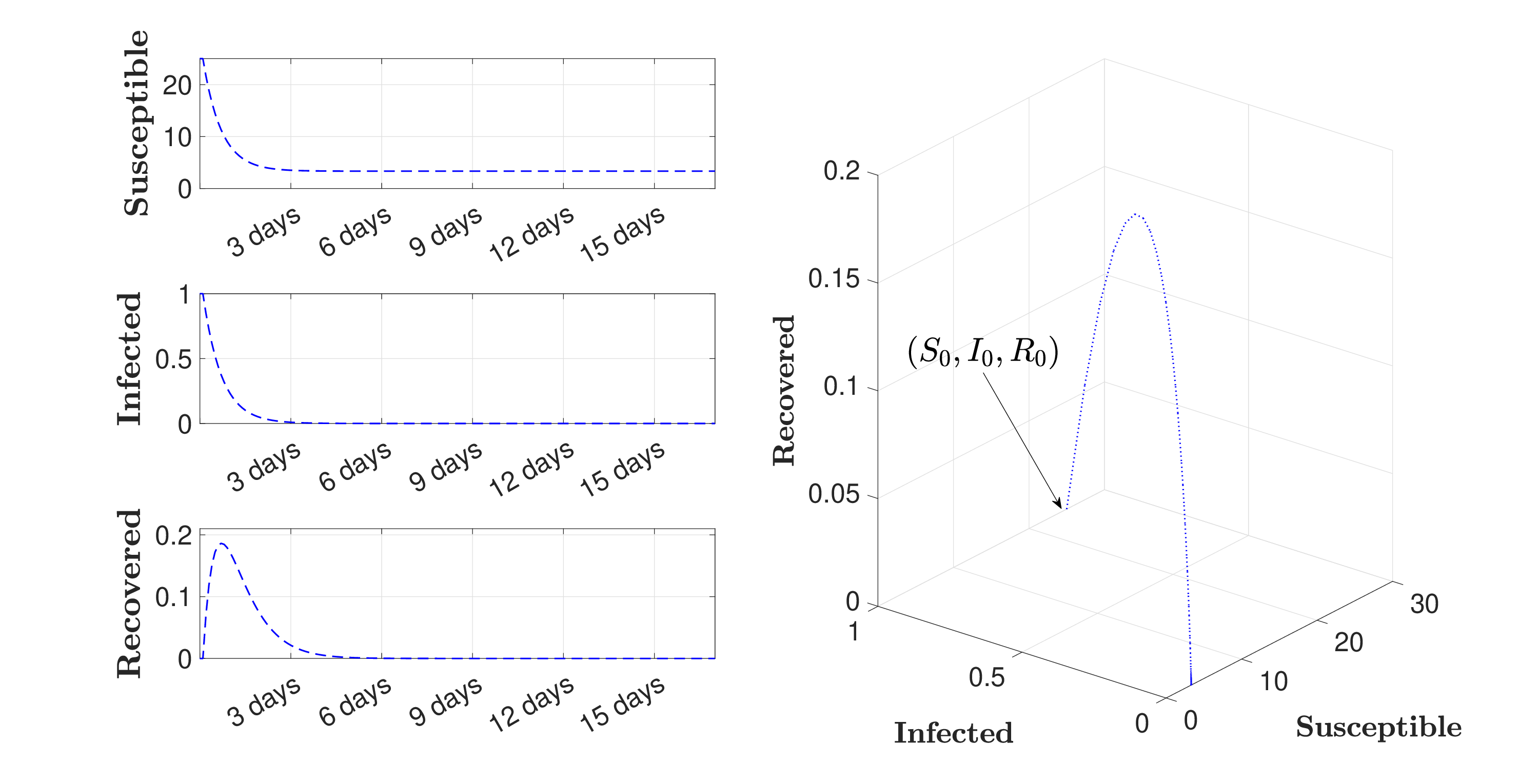}
    \end{center}
    \caption{{\bf The disease is eradicated.} Values of the parameters: $a=1.5$, $b=0.5$, $c=0.7$, $\gamma=1.25$, $q=5$. Initial condition: $(S_0,I_0,R_0)=(25,1,0)$.}
    \label{simulation_theorem_3.1.1}
\end{figure}

In Figure \ref{simulation_theorem_3.1.1}, we take $\gamma=1.25$ and then $\gamma<a+b+c$. Then, the number of infected individuals decreases to zero in approximately 3 days (i.e., the disease is eradicated), the number of susceptible individuals tends to $q/a$ and the number of recovered individuals goes to zero, as we proved in the first statement of Theorem \ref{th-SIR-autonomouos}.

\begin{figure}[H]
    \begin{center}
    \includegraphics[width=\textwidth]{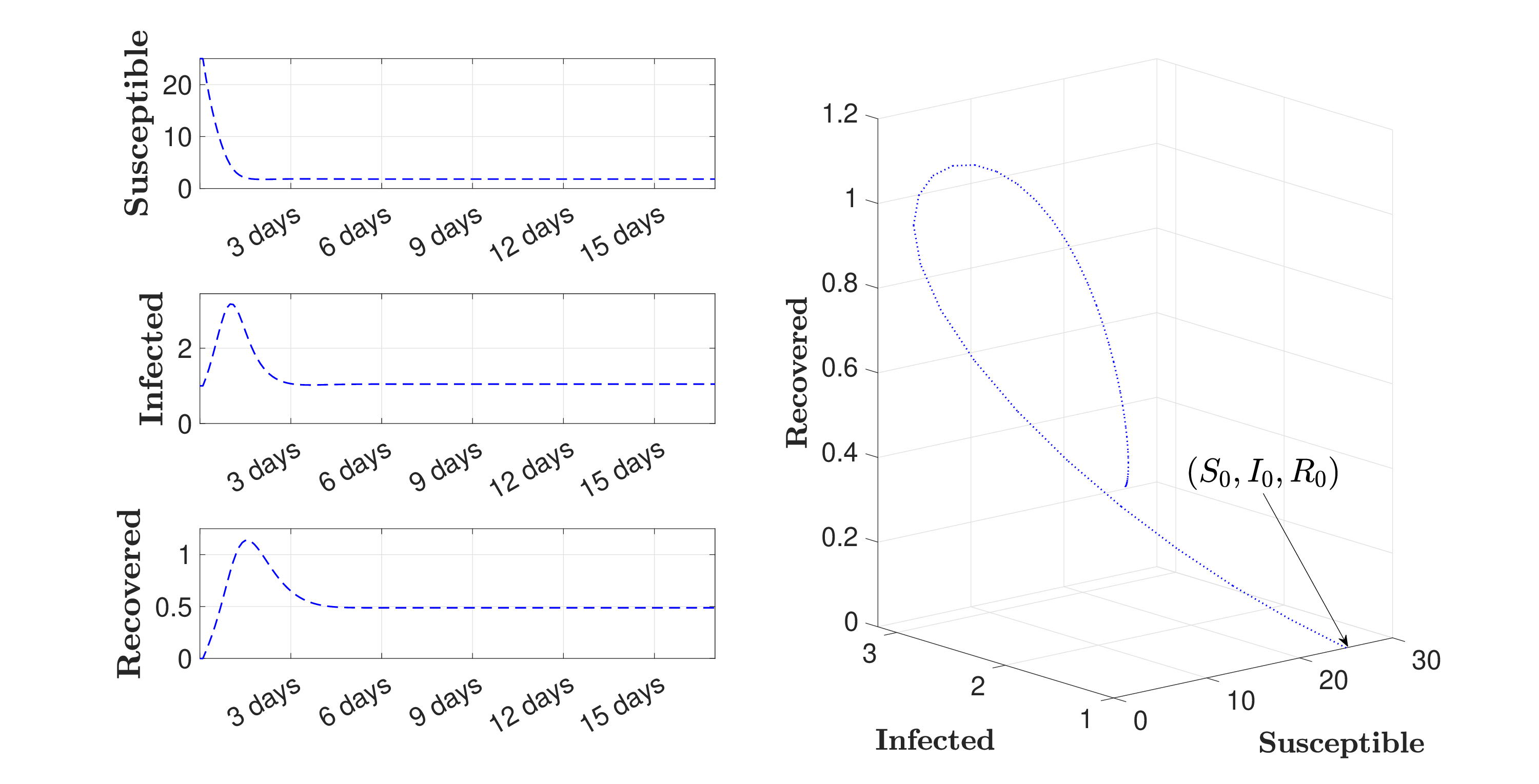}
    \end{center}
    \caption{{\bf The disease becomes endemic}. Values of the parameters: $a=1.5$, $b=0.5$, $c=0.7$, $\gamma=5$, $q=5$. Initial condition: $(S_0,I_0,R_0)=(25,1,0)$.}
    \label{simulation_theorem_3.1.2}
\end{figure}

In Figure \ref{simulation_theorem_3.1.2}, we increase the value of $\gamma=5$. As a consequence, condition $\gamma>a+b+c$ fulfills and therefore the disease becomes endemic, as proved in the second statement of Theorem \ref{th-SIR-autonomouos}. In addition, from the evolution of the infected individuals in the second little panel on the left side of the figure, we can see how the peak of the infection occurs after 1 day and the number of infected individuals is stabilized after 3 days.

\subsection{SIR model with vital dynamics, reinfection and random fluctuations on the diffusion coefficient}\label{sect6.2}

Let us take $d\in(0,\gamma]$ and define $\Phi_\gamma:\Omega\to[-d,d]$ by
\begin{equation}
    \Phi_\gamma(z)=\frac{2d}{\pi}\arctan(z),\label{phi}
\end{equation}
\noindent where $z$ denotes the Ornstein-Uhlenbeck process (see \cite{Caraballo-Kloeden-Schmalfu}).

Then, our goal in this subsection is to depict numerical simulations regarding the random SIR model with reinfection \eqref{odeS}-\eqref{odeR} with $\Phi(\theta_t\omega)=\Phi_\gamma(z(\theta_t\omega))$, which is bounded for every $t\in\R$ and $\omega\in\Omega$ and satisfies condition \eqref{eq-ergodic-condition-for-Phi} in Theorem \ref{th-q-cte-gamma-pertturbed} (see \cite[Proposition 4.1]{Caraballo-Colucci-Cruz-Rapaport-20}). To this end, we set $a=1.5$, $b=0.5$, $c=0.7$, $q=5$, $d=1.5$ and the initial condition $(S_0,I_0,R_0)=(25,1,0)$.

\begin{figure}[H]
    \begin{center}
    \includegraphics[width=\textwidth]{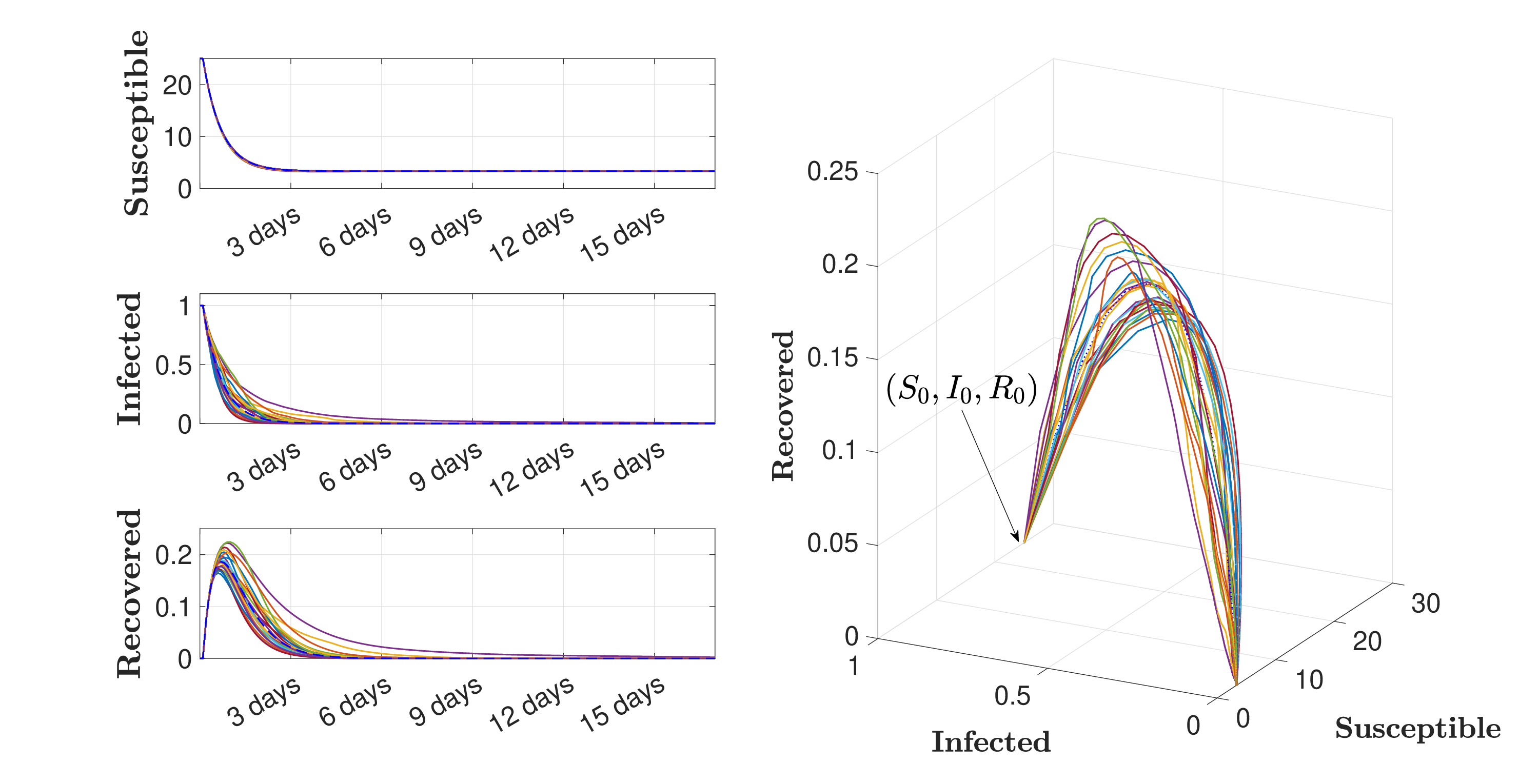}
    \end{center}
    \caption{{\bf The disease is eradicated.} Values of the parameters: $a=1.5$, $b=0.5$, $c=0.7$, $\gamma=1.25$, $q=5$, $d=1.5$. Initial condition: $(S_0,I_0,R_0)=(25,1,0)$.}
    \label{simulation_theorem_4.2}
\end{figure}

In Figure \ref{simulation_theorem_4.2} we choose $\gamma=1.25$, thus $\gamma<a+b+c$, whence the number of infected and recovered individuals decreases to zero (i.e., the disease is eradicated, in less than 3 days) and the number of susceptible individuals tends to $q/a$, as proved in Theorem \ref{th-q-cte-gamma-pertturbed}.

\begin{figure}[H]
    \begin{center}
    \includegraphics[width=\textwidth]{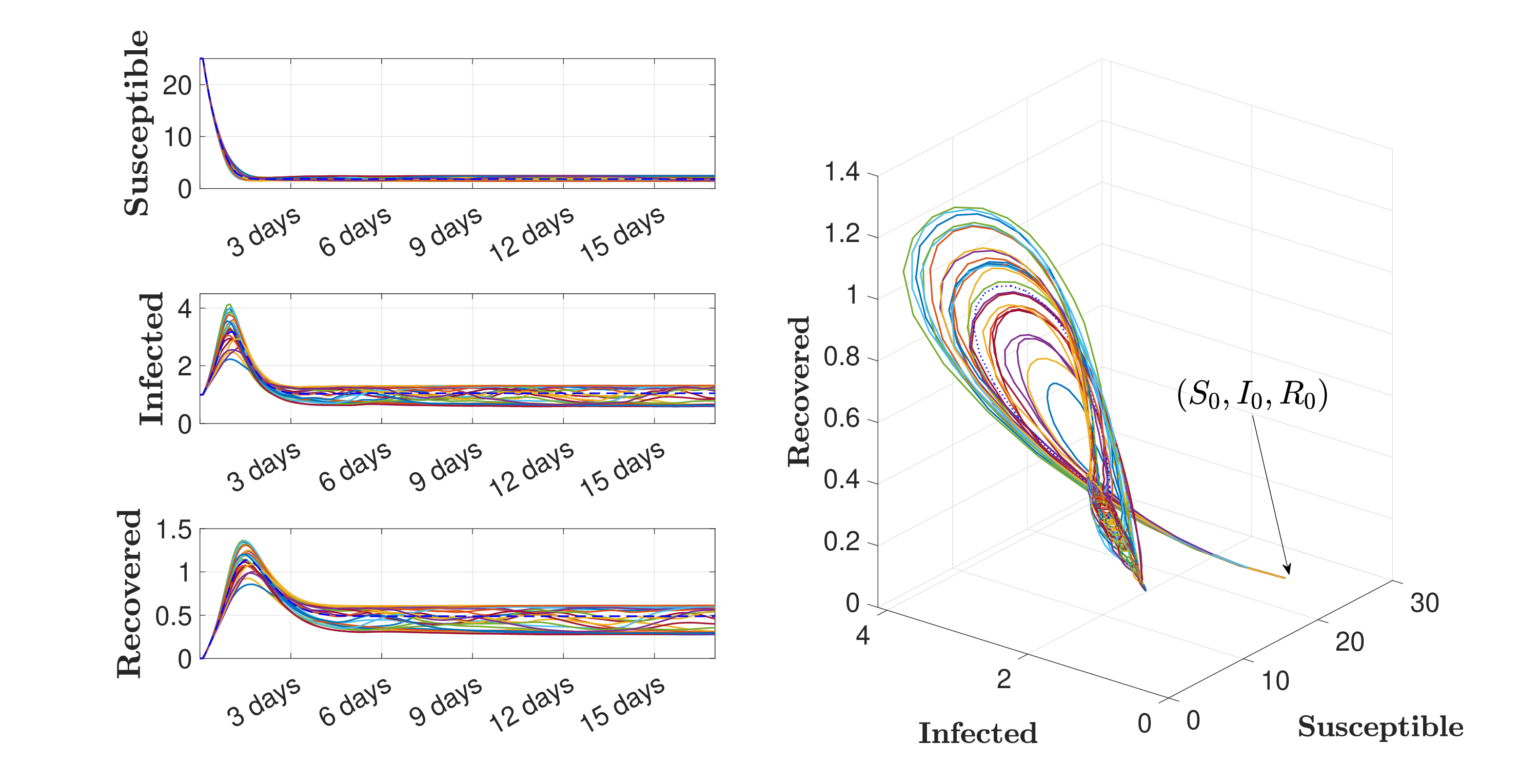}
    \end{center}
    \caption{{\bf The disease becomes endemic}. Values of the parameters: $a=1.5$, $b=0.5$, $c=0.7$, $\gamma=5$, $q=5$, $d=1.5$. Initial condition: $(S_0,I_0,R_0)=(25,1,0)$.}
    \label{simulation_theorem_4.3}
\end{figure}

In Figure \ref{simulation_theorem_4.3} we increase the value of $\gamma=5$. As a result, condition $\gamma_0>a+b+c$ is fulfilled, where $\gamma_0=3.5$ and the disease becomes endemic, as we proved in Theorem \ref{th-condition-below-boundness-gamma}.

\subsection{SIR model with vital dynamics, reinfection and random fluctuations on the
diffusion coefficient and the recruitment rate}\label{section6.3}

Let $d\in(0,\gamma]$ and $e\in(0,q]$ and consider $\Phi_\gamma:\Omega\to[-d,d]$ and $\Phi_q:\Omega\to[-e,e]$ given as follows
\begin{equation}
    \Phi_\gamma(z)=\frac{2d}{\pi}\arctan(z),\quad \Phi_q(z)=\frac{2e}{\pi}\arctan(z)\label{phi1}
\end{equation}
\noindent where $z$ denotes the Ornstein-Uhlenbeck process (see \cite{Caraballo-Kloeden-Schmalfu}).

Then, our goal in this subsection is to depict numerical simulations regarding the random SIR model with reinfection \eqref{eS1}-\eqref{eR1} with $\Phi(\theta_t\omega)=\Phi_\gamma(z(\theta_t\omega))$ and $q(t,\omega)=\Phi_q(z(\theta_t\omega))$, which are bounded for every $t\in\R$ and $\omega\in\Omega$.

We set $a=1.5$, $b=0.5$, $c=0.7$, $q=5$, $d=1.5$, $e=0.5$ and $(S_0,I_0,R_0)=(25,1,0)$.

\begin{figure}[H]
    \begin{center}
    \includegraphics[width=\textwidth]{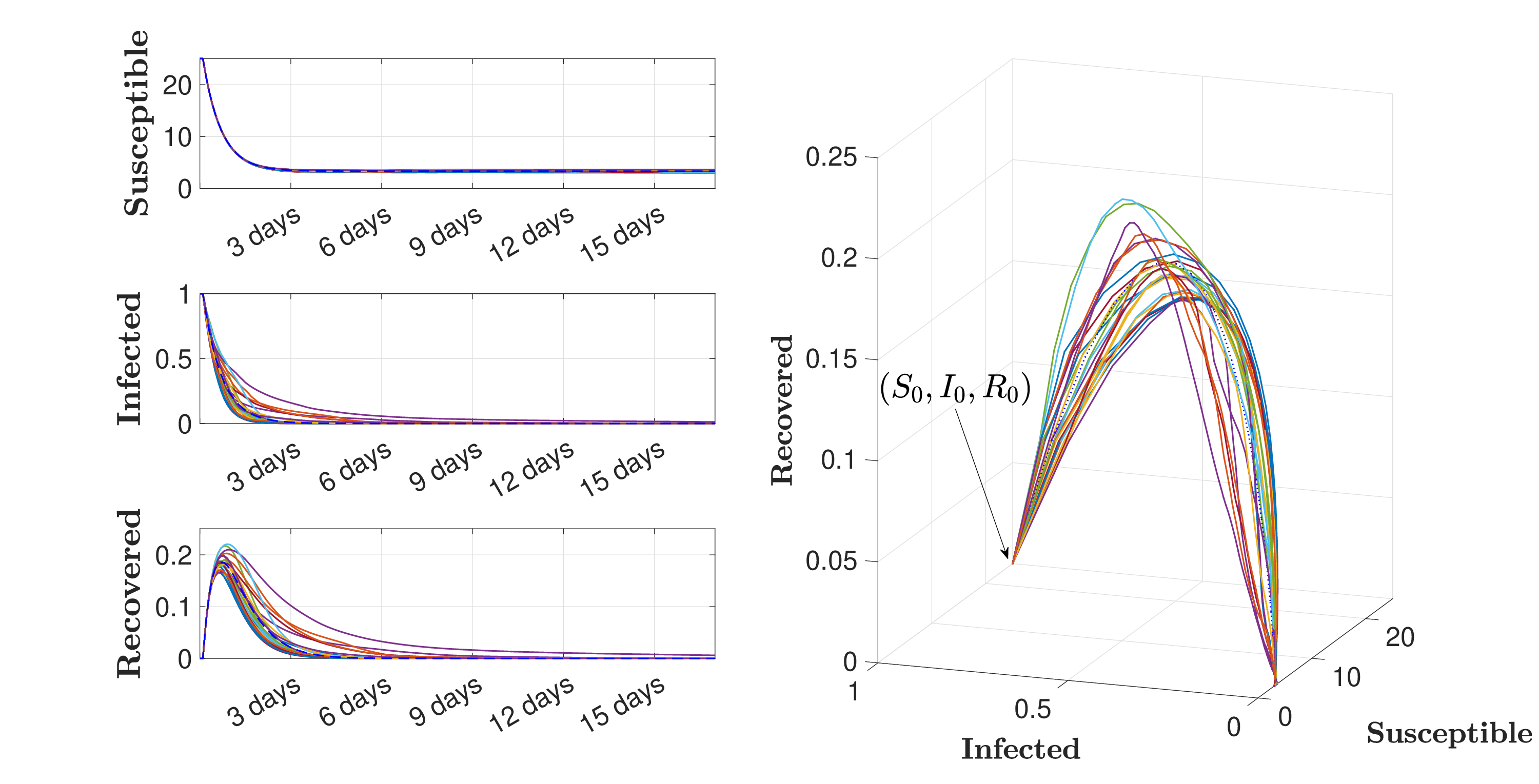}
    \end{center}
    \caption{{\bf The disease is eradicated.} Values of the parameters: $a=1.5$, $b=0.5$, $c=0.7$, $\gamma=1.25$, $q=5$, $d=1.5$, $e=0.5$. Initial condition: $(S_0,I_0,R_0)=(25,1,0)$.}
    \label{simulation_theorem_5.2.1}
\end{figure}

In Figure \ref{simulation_theorem_5.2.1} we consider $\gamma=1.25$, hence condition $\gamma<a+b+c$ holds and the disease is eradicated, in less than 3 days.

\begin{figure}[H]
    \begin{center}
    \includegraphics[width=\textwidth]
    {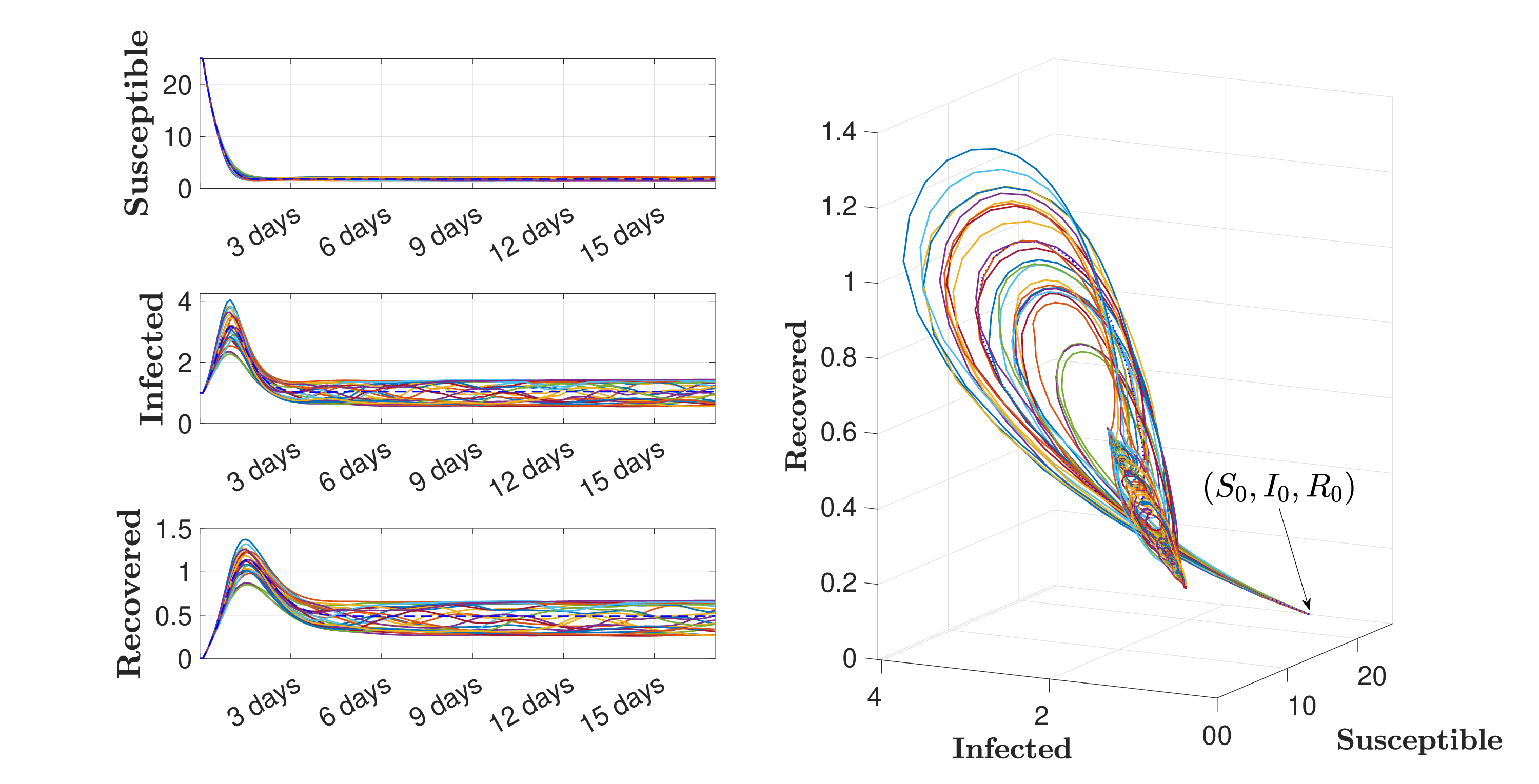}
    \end{center}
    \caption{{\bf The disease becomes endemic}. Values of the parameters: $a=1.5$, $b=0.5$, $c=0.7$, $\gamma=5$, $q=5$, $d=1.5$, $e=0.5$. Initial condition: $(S_0,I_0,R_0)=(25,1,0)$.}
    \label{simulation_theorem_5.2.2}
\end{figure}

Nevertheless, in Figure \ref{simulation_theorem_5.2.2} we increase the value of $\gamma=5$. As a result, $\gamma_0q_0/q_1>a+b+c$. Then, the number of infected individuals remains strictly positive (Theorem \ref{theorem_cases_q_no_constant}).

\section{Conclusion}\label{conclusions}

In this paper, we present the classical SIR model with vital dynamics and reinfection, where we also introduce random and nonautonomous perturbations in an attempt to obtain more realistic systems that describe real epidemics in a better manner.
We study the long-time dynamics of the disease in every resulting system by using the theory of nonautonomous and random dynamical systems. We obtained conditions under which the disease is eradicated or becomes endemic, the existence of attractors, and hyperbolic solutions associated with these models. 

We investigate in detail several aspects of disease modeling within the context of random and nonautonomous systems, motivated by situations in real life. We have proved that the recruitment rate $q$ characterizes the nature of the absorbing sets and the attractors, shedding light on the disease-free solution under different conditions. Apart from that, the stochastic nature of $\gamma$ has a significant influence on the structure of attractors and the hyperbolicity of global solutions. Moreover, we highlighted the importance of considering random perturbations in the diffusion coefficient $\gamma$, reflecting the natural movements of individuals within confined spaces such as schools, supermarkets, and neighborhoods. 


Our exploration has also uncovered novel aspects such as random transmission, reinfection, and alternative modeling techniques. In addition, we have questioned the possibility of defining a proper basic reproductive number in nonautonomous or random frameworks, while also acknowledging specific cases where such definition have been feasible (remarks \ref{remark_3.1} and \ref{Remark_5.1}) or not (Remark \ref{Remark_reproductive_number_fails}). 






\section*{Financial disclosure}

This work has been partially supported by the Spanish Ministerio de Ciencia e Innnovaci\'on, Agencia Estatal de Investigaci\'on (AEI) and Fondo Europeo de Desarrollo Regional (FEDER) under the project PID2021-122991NB-C21, the S\~ao Paulo Research Foundation (FAPESP) under grants 2022/00176-0, 2017/21729-0 and 2018/10633-4 and the CAPES-PRINT under grant 88887.912449/2023-00.


\bibliographystyle{abbrv}
\bibliography{references_SIR}

\begin{thebibliography}{10}

\bibitem{Albani-Zubelli}
V.~V.~L. Albani and J.~P. Zubelli.
\newblock Stochastic transmission in epidemiological models.
\newblock {\em Journal of Mathematical Biology}, 88(25):1432--1416, 2024.

\bibitem{Aragao-Costa-2011}
E.~R. Aragão-Costa, T.~Caraballo, A.~N. Carvalho, and J.~A. Langa.
\newblock Stability of gradient semigroups under perturbations.
\newblock {\em Nonlinearity}, 24(7):2099, jun 2011.

\bibitem{Aragao-Costa-2013}
E.~R. Aragão-Costa, T.~Caraballo, A.~N. Carvalho, and J.~A. Langa.
\newblock Non-autonomous {M}orse-decomposition and {L}yapunov functions for
  gradient-like processes.
\newblock {\em Transactions of the American Mathematical Society},
  365:5277--5312, 04 2013.

\bibitem{Barreira-Valls-upper-tempered}
L.~Barreira and C.~Valls.
\newblock Tempered exponential behavior for a dynamics in upper triangular
  form.
\newblock {\em Eletronic Journal of Qualitative Theory of Differential
  Equations}, 2018(77):1 -- 22, 2018.

\bibitem{Battateli-Palmer-2015}
F.~Battelli and K.~J. Palmer.
\newblock Criteria for exponential dichotomy for triangular systems.
\newblock {\em Journal of Mathematical Analysis and Applications},
  428(1):525--543, 2015.

\bibitem{daniel-bernoulli-1960}
D.~Bernoulli.
\newblock Essai d'une nouvelle analyse de la mortalit\'e cause par la petite
  v\'erole, et des aventages de l'inoculation pour la pr\'evenir.
\newblock {\em M\'emoires de math\'ematique et de physiques tires des registres
  de l'Academie Royale des Sciences}, pages 1--75, 1760.

\bibitem{Bortolan-Carvalho-Langa-book}
M.~C. Bortolan, A.~N. Carvalho, and J.~A. Langa.
\newblock {\em Attractors Under Autonomous and Non-autonomous Perturbations}.
\newblock Mathematical Surveys and Monographs. American Mathematical Society,
  Providence RI, 2020.

\bibitem{Caraballo-Carvalho-Langa-OliveiraSousa-2021}
T.~Caraballo, A.~N. Carvalho, J.~A. Langa, and A.~N. Oliveira-Sousa.
\newblock The effect of a small bounded noise on the hyperbolicity for
  autonomous semilinear differential equations.
\newblock {\em Journal of Mathematical Analysis and Applications}, 500:125134,
  2021.

\bibitem{Caraballo-Colucci-17}
T.~Caraballo and R.~Colucci.
\newblock A comparison between random and stochastic modeling for a {SIR}
  model.
\newblock {\em Communications on Pure and Applied Analysis}, 16(1):151--162,
  2017.

\bibitem{Caraballo-Colucci-Cruz-Rapaport-20}
T.~Caraballo, R.~Colucci, J.~L\'opez de~la Cruz, and A.~Rapaport.
\newblock Study of the chemostat model with non-monotonic growth under random
  disturbances on the removal rate.
\newblock {\em Mathematical Biosciences and Engineering}, 17:7480--7501, 04
  2020.

\bibitem{Caraballo-Kloeden-Schmalfu}
T.~Caraballo, P.~Kloeden, and B.~Schmalfu{\ss}.
\newblock Exponentially stable stationary solutions for stochastic evolution
  equations and their perturbation.
\newblock {\em Applied Mathematics and Optimization}, 50:183--207, 10 2004.

\bibitem{Caraballo-Carvalho-Oliveira-Sousa-NRA}
T.~Caraballo, J.~A. Langa, A.~N. Carvalho, and A.~N. Oliveira-Sousa.
\newblock Continuity and topological structural stability for nonautonomous
  random attractors.
\newblock {\em Stochastics and Dynamics}, 22(07):2240024, 2022.

\bibitem{Carvalho-Langa-Robison-book}
A.~Carvalho, J.~Langa, and J.~Robinson.
\newblock {\em Attractors for infinite-dimensional non-autonomous dynamical
  systems}, volume 182.
\newblock Springer, New York, 2013.

\bibitem{Carvalho-Langa-2}
A.~N. Carvalho and J.~A. Langa.
\newblock Non-autonomous perturbation of autonomous semilinear differential
  equations: Continuity of local stable and unstable manifolds.
\newblock {\em Journal of Differential Equations}, 233(2):622--653, 2007.

\bibitem{Chueshov}
I.~Chueshov.
\newblock {\em Monotone Random Systems theory and applications}.
\newblock Springer-Verlag, New York, 2002.

\bibitem{calculoR0}
P.~V. den Driessche and J.~Watmough.
\newblock Reproduction numbers and sub-threshold endemic equilibria for
  compartmental models of disease transmission.
\newblock {\em Mathematical biosciences}, 180(1):29--48, 2002.

\bibitem{calculoR0old}
O.~Diekmann, J.~Heesterbeek, and J.~Metz.
\newblock On the definition and the computation of the basic reproduction ratio
  $\mathcal{R}_0$ in models for infectious diseases in heterogeneous
  populations.
\newblock {\em Journal of Mathematical Biology}, 28:365--382, 1990.

\bibitem{william-heaton-hamer-1906}
W.~H. Hamer.
\newblock The milroy lectures on epidemic disease in england - the evidence of
  variability and persistence of type.
\newblock {\em The Lancet}, 1(4305):733--739, 1906.

\bibitem{kermack-mckendrick-1927}
W.~O. Kermack and A.~G. McKendrick.
\newblock A contribution to the mathematical theory of epidemics.
\newblock {\em Proceedings of the Royal Society of London}, 115:700--721, 1927.

\bibitem{Kloeden-Kozyakin}
P.~E. Kloeden and V.~S. Kozyakin.
\newblock The dynamics of epidemiological systems with nonautonomous and random
  coefficients.
\newblock {\em MESA}, 2(2):159--172, 2011.

\bibitem{Kloeden-Potzsche}
P.~E. Kloeden and C.~Pötzsche.
\newblock Nonautonomous bifurcation scenarios in {SIR} models.
\newblock {\em Mathematical Methods in the Applied Sciences},
  38(16):3495--3518, 2015.

\bibitem{Martcheva}
M.~Martcheva.
\newblock {\em An Introduction to Mathematical Epidemiology}, volume~61 of {\em
  Texts in Applied Mathematics}.
\newblock Springer New York, NY, 2015.

\bibitem{miranda2022}
A.~Miranda-Mendizabal, S.~Recoder, E.~C. Sebastian, M.~C. Closas, D.~L. Ureña,
  R.~Manolov, N.~M. Santander, C.~G. Forero, and P.~Castellv\'i.
\newblock Socio-economic and psychological impact of covid-19 pandemic in a
  spanish cohort bioval-d-covid-19 study protocol.
\newblock {\em Gaceta Sanitaria}, 36:70--73, 2022.

\bibitem{ronald-ross-1911}
R.~Ross.
\newblock {\em The prevention of Malaria, 2nd edition}.
\newblock John Murray, London, 1911.

\end{thebibliography}
\end{document}